\documentclass[12pt]{amsart}
\usepackage{graphicx,color}
\usepackage{amssymb}
\usepackage{tikz-cd}
\let\temp\emptyset
\usepackage{mathabx}
\let\emptyset\temp
\usepackage{verbatim}
\usepackage[colorlinks=false]{hyperref}
\usepackage[textsize=small]{todonotes}
\usepackage{enumitem}
\usepackage[margin=1.5in]{geometry}

\setlist{leftmargin=2em}

\newtheorem*{mainthm}{Main Theorem}

\newtheorem*{serieslaw}{Series Law}
\newtheorem*{parallellaw}{Parallel Law}

\author {Jeremy Kahn}
\title[bounds for bounded-primitive type]{A priori bounds for some\\
   infinitely renormalizable quadratics:\\
   {\small I. Bounded primitive combinatorics}}
\newcommand{\from}{\mathpunct:}
\numberwithin{equation}{section}
\newcommand{\ceil}[1]{\lceil #1 \rceil}

\newtheorem{thm}{Theorem}[equation]
\newtheorem{cor}[equation]{Corollary}
\newtheorem{lem}[equation]{Lemma}
\newtheorem{prop}[equation]{Proposition}
\newtheorem{theorem}[equation]{Theorem}
\newtheorem{lemma}[equation]{Lemma}

\newtheorem{corollary}[equation]{Corollary}
\theoremstyle{remark}
\newtheorem{rem}[equation]{Remark}

\theoremstyle{definition}
\newtheorem{defn}[equation]{Definition}

\numberwithin{figure}{section}

\def\note#1{}

\newcommand{\di}{\partial}

\newcommand{\ra}{\rightarrow}

\newcommand{\imply}{\Rightarrow}

\def\ssk{\smallskip}
\def\msk{\medskip}

\def\sm{\smallsetminus}

\newcommand{\arrow}{\longrightarrow}

\newcommand{\dist}{\operatorname{dist}}

\newcommand{\cl}{\operatorname{cl}}
\newcommand{\inter}{\operatorname{int}}
\renewcommand{\mod}{\operatorname{mod}}

\newcommand{\tl}{\tilde}

\newcommand{\can}{{\mathrm{can}}}

\newcommand{\dominate}{\multimap}

\newcommand{\ta}{\leadsto}

\newcommand{\supp}{\operatorname{supp}}

\renewcommand{\d}{{\diamond}}

\newcommand{\eps}{{\varepsilon}}

\newcommand{\De}{{\Delta}}
\newcommand{\de}{{\delta}}
\newcommand{\la}{{\lambda}}

\newcommand{\si}{{\sigma}}

\newcommand{\om}{{\omega}}

\newcommand{\bal}{{\mbox{\boldmath$\alpha$} }}
\newcommand{\ba}{{\bold {a}  }}

\newcommand{\bbe}{{\mbox{\boldmath$\beta$} }}

\newcommand{\AAA}{{\mathcal A}}

\newcommand{\DD}{{\mathcal D}}
\newcommand{\EE}{{\mathcal E}}

\newcommand{\II}{{\mathcal I}}
\newcommand{\FF}{{\mathcal F}}
\newcommand{\GG}{{\mathcal G}}

\newcommand{\HH}{{\mathcal H}}
\newcommand{\KK}{{\mathcal K}}
\newcommand{\LL}{{\mathcal L}}

\newcommand{\OO}{{\mathcal O}}
\newcommand{\PP}{{\mathcal P}}

\newcommand{\WW}{{\mathcal W}}

\newcommand{\A}{{\Bbb A}}
\newcommand{\C}{{\Bbb C}}

\newcommand{\D}{{\Bbb D}}

\newcommand{\M}{{\Bbb M}}
\newcommand{\N}{{\Bbb N}}

\newcommand{\R}{{\Bbb R}}
\newcommand{\T}{{\Bbb T}}

\newcommand{\Z}{{\Bbb Z}}

\newcommand{\hor}{{\mathrm{h}}}
\newcommand{\ver}{{\mathrm{v}}}

\newcommand{\g}{{\bf g}}
\newcommand{\h}{{\bf h}}

\def\BG{{\mathbf{G}}}

\def\BH{{\mathbf{H}}}

\def\B0{{\mathbf{0}}}

\def\BQ{{\mathbf{Q}}}

\def\BV{{\mathbf{V}}}
\def\BU{{\mathbf{U}}}

\newcommand{\area}{\operatorname{area}}
\newcommand{\Flag}{\operatorname{Flag}}

\renewcommand{\max}{{\mathrm{max}}}

\newcommand{\fjset}{FJ-set}

\catcode`\@=12

\def\Empty{}
\newcommand\oplabel[1]{
  \def\OpArg{#1} \ifx \OpArg\Empty {} \else
  	\label{#1}
  \fi}

\def\makeabbrevs{\def\o{\omega}\def\g{\gamma}\def\G{\Gamma}\def\h{\hat}\def\d{\delta}\def\D{\Delta}\def\O{\Omega}\def\b{\beta}\def\l{\lambda}}

\renewcommand{\mod}{\operatorname{mod}}
\newcommand{\phl}[1]{\left\vert #1 \right\vert}

\newcommand{\wcan}{W_{\operatorname{can}}}
\let\Wcan\wcan
\newcommand{\thh}{^{\text{th}}}
\newcommand{\Sym}{\operatorname{Sym}}

\def\Int{\operatorname{Int}}
\def\Cl{\operatorname{Cl}}
\def\from{\mathpunct:}
\def\AA{{\mathcal A}}

\def\ezero{\epsilon_0}
\def\M{\core M}
\def\core#1{{\boldsymbol{#1}}}
\def\bfalpha{\core\alpha}
\newcommand{\Hull}{\operatorname{Hull}}

\def\pairing#1{\left< #1 \right>}
\newcommand{\pair}[2]{\left< #1, #2 \right>}
\newcommand{\abs}[1]{\left\vert #1 \right\vert}

\newcommand{\norm}[1]{\left\Vert #1 \right\Vert}
\newcommand{\normone}[1]{\norm{#1}_1}

\newcommand{\ff}{\mathbf f}
\newcommand{\bgg}{\mathbf g}
\newcommand{\dout}{\partial_{\text{out}}}
\renewcommand{\cl}[1]{\overline{#1}}
\newcommand{\inv}[1]{^{-1}(#1)}
\newcommand{\nnm}{\dotdiv}
\newcommand{\chat}{\hat \C}

\newcommand{\hx}[1]{X^{#1}}
\newcommand{\doms}{\multimap}

\renewcommand{\ba}{\boldsymbol{\alpha}}

\def\norm#1{\Vert #1 \Vert}
\newcommand{\lollypop}{\multimap}
\newcommand{\poincare}{Poincar\'e}

\newcommand{\pg}{\pi_G}
\newcommand{\ph}{\pi_H}

\thanks{
   This work was supported in part by Sloan Research Fellowship
 and NSF grants DMS-8920768 and DMS-9022140.}
\date{\today}

\begin{document}

\maketitle

\begin{abstract}
	We prove the \emph{a priori} bounds for infinitely renormalizable quadratic polynomials for which we can find an infinite sequence of primitive renormalizations such that the ratios of the periods of successive renormalizations is bounded.
\end{abstract}

\tableofcontents

\newcommand{\smallheading}[1]{\msk {\it #1.}}
\renewcommand{\smallheading}[1]{\subsubsection{#1}}

\section{Introduction} \label{sec:intro}
We will introduce the Main Theorem in Section \ref{intro:main}
(and list its many consequences in Section \ref{intro:consequences}), after fixing notation in Section \ref{intro:defs}.
We then outline the proof in Section \ref{intro:outline}, 
and finish this section with some general notation and the acknowledgements. 
\subsection{Definitions and notation for renormalization} \label{intro:defs}
We review some basic definitions in order to fix notation throughout the paper. 
With one minor but significant exception\footnote{We refer to the quadratic-like map $f^p\from U \to V$ as a ``$p$-renormalization of $f$''
rather than a ``renormalization of $f^p$''.}, 
the notation follows \cite{McM},
and the reader is referred to this book for an excellent introduction to complex dynamics and renormalization.

A  {\it topological disk} (or simply ``disk'') means a simply connected  domain in some Riemann surface $S$.
A \emph{polynomial-like map} is a proper holomorphic map $g\from U \to V$ where $U \Subset V \subset \C$.
A quadratic-like map is a polynomial-like map of degree 2.
We denote the filled-Julia set of $g$ by $K(g)$, or just $K$ when the context is clear. 

A quadratic polynomial  $f\from z\mapsto z^2+c$
is called {\it renormalizable with period $p$}
(or \emph{$p$-renormalizable})
if there exist disks $U, V \in \C$ with $0 \in U$ such that $f^p\from U\to V$
is a quadratic-like map with connected Julia set.
This quadratic-like map is called a\footnote{We use `a' rather than `the' because the domain and range of the renormalization is not uniquely determined. 
On the other hand, the germ of the renormalization at its filed-Julia set is uniquely determined, 
so we (and others) will sometimes use `the' when the focus is the germ or the small Julia set. } 
{\it renormalization} of $f$ (of period $p$).
While $U$ and $V$ are not uniquely determined,
the filled Julia set $K_p$ of the $p$-renormalization is determined\footnote{See Theorem 7.1 of \cite{McM}} 
by $f$ and $p$. 
As in \cite{McM},
we let $K_p(i) = f^i(K_p)$ for $i = 0, \ldots, p-1$, 
and we let $\KK_p = \bigcup_i K_p(i)$. 
In a context where $p$ is clear we will sometimes write $\KK_p$ as $\KK$.
We call $K_p$ the \emph{central small Julia sets};
we call the $K_p(i)$ the small Julia sets,
and $\KK_p$ the union of small Julia sets;
although these are actually \emph{filled Julia sets},
we will allow ourselves this small elision. 

We say that the $p$-renormalization of $f$ is \emph{primitive} if the small Julia sets of period $p$ are pairwise disjoint.
In this case,
we let $\gamma_p \equiv \gamma_p(f)$ denote the closed geodesic around $K_p$ in $\C \sm \KK_p$ in the hyperbolic (\poincare) metric,
and likewise let $\gamma_p(i)$ denote the geodesic around $K_p(i)$. 
If $g\from U\to V$ is quadratic-like,
we let $\gamma(g)$ (or just $\gamma$) be the core geodesic of the annulus $V \sm K(g)$. 
If $\eta$ is any closed geodesic,
we let $\phl \eta$ denote the length of $\eta$.

If there is an infinite
sequence of periods $p_0< p_1<\dots$ such that $f$ is primitively
$p_k$-renormalizable then  $f$ is called {\it infinitely primitively renormalizable}.
If additionally, there exists a $B$ such that $p_{k+1} / p_k\leq B$, $k=0,1,\dots$,
then  $f$ is called infinitely renormalizable of $B$-bounded primitive type.
Such a map has \emph{a priori} bounds if there exists an $\eps>0$ and
a sequence of quadratic-like renormalizations $f^{p_k}\from U_k\to V_k$ such that
$\mod(V_k\sm U_k)\geq \eps$.

In situations where we are interested in the sequence of renormalizations of periods $p_1, p_2, \ldots$,
we let $K_{(n)}, \KK_{(n)}, \gamma_{(n)}$ be shorthand for $K_{p_n}, \KK_{p_n}, \gamma_{p_n}$.
The reader should pay close attention\footnote{Usually the period will be denoted by $p$ (or $q$), while the index will be $k$ or $n$, so one does not have to pay \emph{very} close attention.},
because $K_{n}$ means the small Julia set for the period $n$ renormalization,
while $K_{(n)}$ means the same for the $\underline{n\thh}$ renormalization.

\subsection{Statement of the main result} \label{intro:main}
We can now state the main result of this paper.
\begin{mainthm}
Let $f$ be an infinitely renormalizable
quadratic polynomial of bounded primitive type. Then $f$ has a priori bounds.
\end{mainthm}

In subsequent papers \cite{decorations, molecules}
the author and M. Lyubich proved a priori bounds for 
a  class of infinitely renormalizable maps of unbounded type.
In very recent work \cite{satellite}, D. Dudko and M. Lyubich,
building on the work of this paper,
remove the hypothesis of primitive renormalization, 
to the conclude that the Main Theorem holds for all infinitely renormalizable quadratic polynomials of bounded type.

For real quadratics of bounded type \emph{a priori} bounds were proved by Sullivan \cite{S},
see also \cite{LS,LY,MS}.
They were also proved for a class of complex combinatorics of ``high bounded type'' \cite{puzzle}.

\subsection{Consequences} \label{intro:consequences}

The \emph{a priori} bounds have numerous consequences. Let us list some of them
(below $f_c\from c\mapsto z^2+c$ stands for an infinitely primitively renormalizable
quadratic polynomial of bounded type):
\begin{itemize}
\item
{\it The Mandelbrot set is locally connected at $c$, or equivalently,
the polynomial $f_c$ is combinatorially rigid} (see \cite{puzzle}).
The conjecture of local connectivity of the Mandelbrot set (the {\it MLC Conjecture})
formulated  about 40 years ago by Douady and Hubbard  (see \cite{Orsay})
is a central open problem in holomorphic dynamics. Previously, it was established
for all quadratic maps which are not infinitely renormalizable (Yoccoz, see \cite{H}) and
for the class  of infinitely renormalizable maps of high type mentioned above
(see \cite{puzzle}).
\item
   {\it The Julia set $J(f_c)$ is locally connected} (see \cite{HJ,J}).
\item
  {\it  The Feigenbaum-Coullet-Tresser Renormalization Conjecture is valid for primitive combinatorics}.
This conjecture was established in the work of Sullivan \cite{S}, McMullen \cite{towers}
and Lyubich \cite{universe} assuming {\it a priori} bounds (and thus, unconditionally,
for real maps). Now, these results become unconditional for arbitrary primitive combinatorics.
\item
 {\it  Universality and Hairiness of the Mandelbrot set at} $c$.
These properties were conjectured by Milnor \cite{M} and proved in \cite{universe} for maps
with {\it a priori} bounds.
\item
 {\it The Basic Trichotomy for the measure and Hausdorff dimension of the Julia set $J(f_c)$} which was established in \cite{AL}
for maps with {\it a priori} bounds.
\end{itemize}

\subsection{Outline of the proof} \label{intro:outline}
We will now give a brief top-down outline of the proof of the Main Theorem.
We remind the reader that we have introduced the essential notation \emph{once and for all} in Section \ref{intro:defs}.

\smallheading{General strategy: induction on the lengths of the hyperbolic geodesics} It has long been known\footnote{And first appeared in print as Theorem 4.10 in \cite{towers}}
that the \emph{a priori} bounds are equivalent to the assertion that the lengths of the $\gamma_{(n)}$
are bounded. 
Our strategy towards this end is to show that {\it if the length of some $\gamma_{(n)}$ gets long
then it was even longer before}:
There exist  $M>0$ and  $l>0$
such that
\begin{equation} \label{eq:bad-now-worse-earlier}
  |\gamma_{(n)}| > M \imply |\gamma_{(n-l)}| > 2M.
\end{equation}

\smallheading{Localization ``in principle''}
The immediate and obvious difficulty in proving \eqref{eq:bad-now-worse-earlier} is that the number of components of $\KK_{(n)}$ grows exponentially with $n$, 
and of course is unbounded as $n$ goes to infinity. 
While there has been some success in handling an unbounded number of components,
there is certainly a vast advantage to working with a bounded number,
and this paper is able to make a major step forward with the help of such a localization.

The central idea is that, while $n$ is unbounded, $l$ is bounded;
because we assume $p_{k+1}/p_k \le B$ (for all $k$),
the number of components of $\KK_{(n)}$ that are surrounded by $\gamma_{n-l}$ is then at most $B^l$.
Moreover,
if we let $f_{(n-l)}\equiv f^{p_{n-l}}\from U_{(n-l)} \to V_{(n-l)}$ be a renormalization of period $p_{n-l}$ of $f$,
then this $f_{(n-l)}$ is renormalizable of period $p_n/p_{n-l}$,
and the small Julia sets for this renormalization are exactly the components of $\KK_{(n)}$ surrounded by $\gamma_{(n-l)}$.
\begin{figure}
\includegraphics[width=\textwidth]{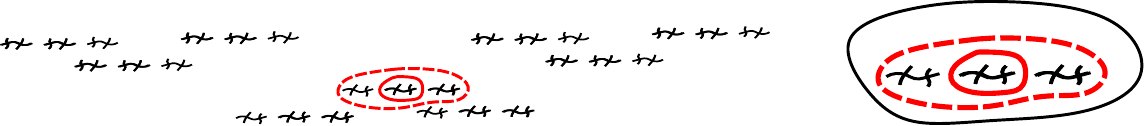}
\caption{On the left we imagine $f$ is renormalizable of periods $3^k$, for $k = 1, 2, 3$. The solid red geodesic is $\gamma_{(3)}$ 
while the dotted red geodesic is shorter than $\gamma_{(2)}$. On the right we have formed a domain $\BV_{(2)}$ on which there is a conformal dynamical system that is renormalizable just of period 3 (the ratio of the last two periods). This dynamical system, formed by the canonical renormalization, has three small Julia sets and two geodesics with the same length as their analogs on the left. To first approximation we can think of it as a quadratic-like renormalization $f_{(2)}$ of $f$.}
\label{fig:local}
\end{figure}
With this in mind, 
we will first state and prove a theorem about a primitively renormalizable quadratic-like map that will be closely related to \eqref{eq:bad-now-worse-earlier}.
\begin{theorem} \label{thm:length-comparison}
There exists $\epsilon > 0$, such that for all $p$ there exists $M$:
Suppose that $g\from U \to V$ is a quadratic-like map that is primitively $p$-renormalizable.
Then (with $\gamma \equiv \gamma(g)$ and $\gamma_p = \gamma_p(g)$), 
\begin{equation} \label{eq:length-comparison}
\phl{\gamma_p}> M \implies \phl{\gamma}\ge \epsilon p \phl{\gamma_p}.
\end{equation}
\end{theorem}

Now suppose $g$ in Theorem \ref{thm:length-comparison} is the renormalization $f_{(n-l)}$ described above,
and $V_{(n-l)}$ is taken to be ``as large as can be'', 
and let $q = p_n/p_{n-l}$. 
Then the $\gamma_q(g)$ in Theorem \ref{thm:length-comparison} is a ``good approximation'' of the $\gamma_{(n)}$ for the original $f$,
and the $\gamma$ for  the $g$ in Theorem \ref{thm:length-comparison} is a ``good approximation'' of the $\gamma_{n-l}$ for the original $f$.
We can choose $q$ \emph{large but bounded},
and \eqref{eq:bad-now-worse-earlier} would follow from \eqref{eq:length-comparison} if these ``good approximations'' were actually equalities.

To make them equalities we introduce (in Section \ref{pseudo-puzzle}) the \emph{pseudo-quadratic-like map} and the \emph{canonical renormalization}.
A pseudo-quadratic-like map $(i, f)\from 
\BU \to \BV$ 
is a generalization of a quadratic-like map 
where the inclusion from $\BU$ to $\BV$ is replaced by an immersion $i$ that has certain properties that imply the (multivalued) $f \circ i^{-1}$ has a (single-valued) quadratic-like restriction. 
We can then,
given our original quadratic polynomial $f$,
define, 
for each $n$,
the \emph{canonical renormalization} $f_{(n)}\from \BU_{(n)}\to \BV_{(n)}$
(really $(i_{(n)}, f_{(n)})$ but we will suppress the $i_{(n)}$ here),
which will be a pseudo-quadratic-like map.
Then,
given $n$ and $l < n$,
and letting $q = p_n/p_{n-l}$,
we'll see that $f_{(n-l)}$ has two crucial properties:
\begin{enumerate}
\item
The geodesic $\gamma(f_{(n-l)})$ (around $K(f_{(n-l)})$ in $\BV_{(n-l)} \sm K(f_{(n-l)})$) has length \emph{equal to} $\gamma_{(n-l)}$.
\item
$f_{(n-l)}$ is renormalizable of period $q$,
and the geodesic $\gamma_q(f_{(n-l)})$ 
(around $K_q(f_{(n-l)})$ in $\BV_{(n-l)} \sm \KK_q(f_{(n-l)})$) has length \emph{greater than or equal to} the length of $\gamma_{(n)}(f)$.  
\end{enumerate}
See Figure \ref{fig:local} for a pictorial representation of these properties.
In Section \ref{sec:conclusions},
we then extend Theorem \ref{thm:length-comparison} to pseudo-quadratic-like maps,
and \eqref{eq:bad-now-worse-earlier} and our Main Theorem then follow.

\begin{figure} 
\includegraphics[width=\textwidth]{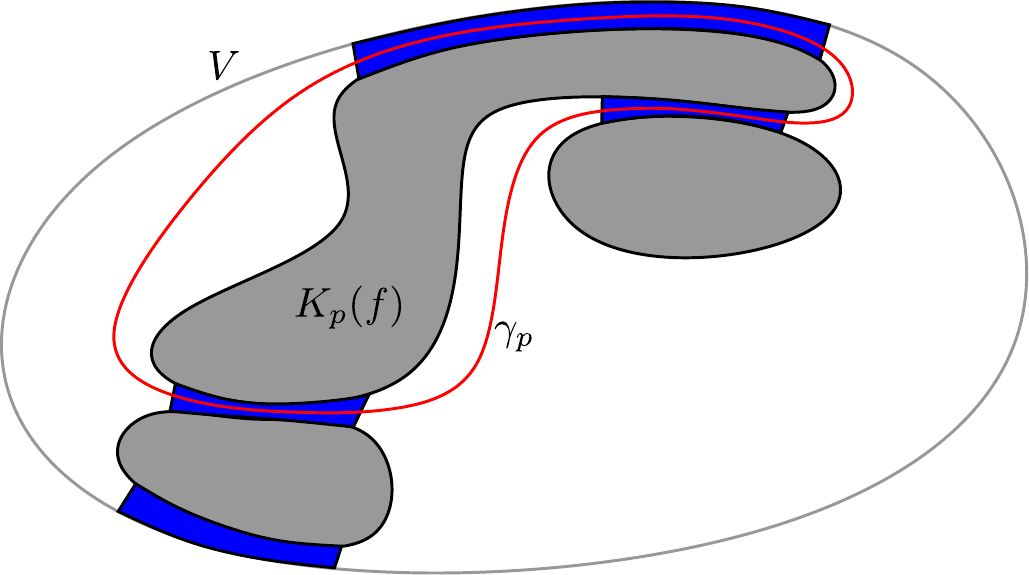}
\caption{Here the gray blobs represent $\KK_p(f)$. The length of the red geodesic $\gamma_p$ going around $K_p(f)$ is roughly the sum of the conformal widths of the thin rectangles, shown in blue.}
\label{fig:thin}
\end{figure}

\smallheading{Canonical weighted arc diagrams} 
Let us now briefly discuss the proof of Theorem \ref{thm:length-comparison}.
In Section \ref{sec:canonical} we develop the theory of the \emph{canonical weighted arc-diagram},
which describes the ``thin rectangles'' of large conformal modulus\footnote{Here we uniformize each of the blue rectangles shown in Figure \ref{fig:thin},
so that the short side goes to a side of length 1. The length of the image of the long side is the conformal modulus.}
(or \emph{width})
which we will then use to estimate $|\gamma_p|$ as the sum of conformal moduli of thin rectangles in $V \sm \KK_p$ crossed by $\gamma_p$. 
(See Figure \ref{fig:thin}).

Roughly speaking, 
the canonical weighted-arc diagram is
$$
     W_\can = \sum W(\alpha)\alpha,
$$
where the $\alpha$ are arcs\footnote{An {\it arc} is a non-trivial homotopy class of properly embedded paths.}
in $V\sm \KK_p$ represented by paths crossing the thin rectangles and the $W(\alpha)$ are their widths. 
Holomorphic maps between Riemann surfaces,
such as covering and inclusions,
then induce relations between their canonical WAD's;
these relations are described in Section \ref{sec:canonical}. 
The most subtle of these relations is 
\emph{domination}, which describes the relation between the canonical WAD's between Riemann surfaces $U \subset V$,
and encodes the parallel and series laws for extremal length.

Theorem \ref{thm:length-comparison} is then reduced (in Section \ref{subsec:hv}) to a statement about the relative weight of \emph{horizontal} and \emph{vertical} arcs in the canonical WAD on $V \sm \KK$, where the horizontal arcs are the ones that go from $\KK$ to $\KK$, and the vertical ones connect $\KK$ to $\partial V$. The idea is that the horizontal arcs contribute only to $|\gamma_p|$ while the vertical arcs contribute to $|\gamma|$; hence Theorem \ref{thm:length-comparison} comes down to showing that the total weight of the horizontal arcs is a definite proportion of the total vertical weight. 

In the broadest terms, 
this relationship between horizontal and vertical is shown using the combinatorial properties of the canonical WAD. 
A much better and more detailed outline is given in Section \ref{subsec:hv-proof-outline};
the reader may wish to skip to that subsection after skimming Section \ref{sec:canonical} and the preliminaries in Section \ref{sec:reduce-to-wcan}.
The relationship itself is proven in Section \ref{sec:hor-vert},
using a combinatorial theorem that is stated and proven in Section \ref{sec:entropy},
after a discussion of Hubbard trees in Section \ref{Hubbard trees}.

\subsection{Terminology and Notation}
We let:\\
$\N =\{0, 1,2,\dots\}$ be the set of natural numbers;\\
$\D=\{z:\,  |z|<1\}$ be the unit disk, $\T_r=\{z: |z|=r\}$, $\T=\T_1$;\\
$\A(r,R)= \{z: r<|z|<R\}$.

We will say a subset $K$ of $\C$ is an \fjset\ (``filled Julia set'')
if $K$ is compact, connected, full, and has more than one element.

\subsection{Acknowledgments}
The original version of this paper was written in collaboration with Mikhail Lyubich following
the mathematical work of the author. The primary result (in particular, Theorem 9.1) of this paper
was then combined with the main result of two subsequent papers 
\cite{decorations,molecules}
of Mikhail
Lyubich and the author to form one result that was the goal of the series
of papers of which this paper was (originally) the first.

This work was partially supported by the Clay Mathematics Institute and the Simons Foundation.
Part of it was done during the author's visit to the IMS at Stony Brook and the Fields Institute in Toronto.
The author is thankful to these Institutions and Foundations.

\section{The canonical weighted arc diagram} \label{sec:canonical}
In this section we introduce the canonical weighted arc-diagram for an open Riemann surface $S$ and describe a number of properties of it. 
This combinatorial object will summarize how close the components of the ideal boundary of $S$ come together;
and the properties will allow us to understand the nearly degenerate parts of the geometry of any complex dynamical system (of bounded topology)  
in terms of combinatorial relations.

\subsection{Proper paths and proper homotopies} \label{subsec:paths}
We let $S$ be a Riemann surface\footnote{So, by definition, $S$ is connected. Up to Section \ref{subsubsec:valid} we use only the smooth (in fact just topological) structure of $S$.}, 
and observe that the following are equivalent:
\begin{enumerate}
\item
$S$ has finitely generated homology,
\item
$S$ has finite genus and finitely many ends,
\item
for some $n$ and $g$, $S$ is diffeomorphic to the surface of genus $g$ with $n$ points (or disjoint closed disks) removed. 
\end{enumerate}
If any (and therefore all) of these hold, we say that $S$ has finite topology.
For the rest of this paper, we will assume that this is the case.
Unless otherwise stated, we will assume that the Euler characteristic $\chi(S)$ is negative.

A (proper) \emph{path} in $S$ is a proper continuous map $a \from I \to S$, where $I \subset \R$ is an open interval. 
Unless otherwise stated, we will always assume (and require) that a path is proper. 
Of course, if there are any paths in $S$, then $S$ is non-compact.
We can \emph{reparametrize} a path $a\from I \to S$ by a homeomorphism $h\from J \to I$, 
and then $a \circ h$ is a \emph{reparametrization} of $a$. 
(Unless otherwise specified, our paths are unoriented, and the homeomorphism $h$ can be orientation-reversing.)
We say that paths $a_0, a_1\from I \to S$ are \emph{homotopic} 
if there exists a proper continuous map $a\from  [0, 1] \times I \to S$ such that $a_i(t) = a(i, t)$ for $i = 0, 1$ (and all $t \in I$).
We say that a path $a\from I \to S$ is \emph{trivial} if there exists a proper continuous map $\hat a\from  [0, 1) \times I \to S$ such that 
$\hat a(0, t) = a(t)$ for all $t \in I$. 
We observe that if a path is trivial, then all homotopic paths are trivial as well.

Let $\EE(S)$ denote the ends of $S$. 
Suppose $a \from I \to S$ is a path, and for simplicity of notation assume that $I \subset \R$ is bounded. 
We then have a unique continuous extension of $a$ to $a\from \cl I \to S \cup \EE(S)$, 
where $a(\di I) \subset \EE(S)$. 
We call $a(\di I)$ the \emph{endpoints} of $a$. 
From this point of view the word ``endpoints'' is quite natural, 
because each endpoint is a point in $\EE(S)$. 
If $S \subset \chat$,
we can identify each end of $S$ with a component of $\chat \sm S$,
and call the appropriate such components the endpoints of $a$. 
From this perspective it may seem somewhat misleading to use the term endpoint, because now each endpoint is a whole component of $\chat \sm S$; for example, it may be a small Julia set.
But we've already chosen the term endpoint, so we will stick with it. 

We will sometimes, by a mild abuse of notation, identify a path with its image.
For example,
we say that two paths are disjoint if their images are disjoint, and that a path is disjoint from some $A \subset S$ if its image is disjoint from $A$.

A \emph{spine} for $S$ is a finite embedded 1-complex $A$ for which there is a deformation retract from $S$ to $A$.
We observe that every finite-topology Riemann surface has a spine.
Let $A$ be any spine for $S$; we observe that the following are equivalent,
for any path $a$ (in $S$):
\begin{enumerate}
\item
$a$ is trivial,
\item
$a$ is homotopic to a path disjoint from $A$,
\item
For any compact set $E \subset S$,
$a$ is homotopic to a path disjoint from from $E$.
\end{enumerate}

We also observe the following, which should alleviate any concerns about ``wild'' paths:
\begin{lemma} \label{lem:un-wild}
Suppose that $A$ is a spine for $S$, and $a_0, a_1\from I \to S$ are paths,
and $I' \subset I$, such that $a_0|_{I'} = a_1|_{I'}$ and $a_i^{-1}(A) \subset I'$ for $i = 0, 1$.
Then $a_0$ and $a_1$ are homotopic.
\end{lemma}
\begin{proof}
First suppose that $I'$ is an initial segment of $I$, and let $I''$ be the complement of $I'$ in $I$.
Then $a_i|_{I'}$ lies in an annular neighborhood of an end of $S$;
we can coordinatize this neighborhood as $[0, \infty) \times \R/\Z$,
and then lift the $a_i|_{I'}$ to $[0, \infty) \times \R$, with the second coordinate going to infinity as we approach the far endpoint of $I''$.
We then let $a_s(t) = (1-s) a_0(t) + s a_1(t)$ in these coordinates,
and push this back down to the annular neighborhood, to obtain a proper homotopy. 
\end{proof}

A more careful argument shows that homotopic embedded paths are isotopic through embedded paths, using the same statement for paths on a compact oriented surface with boundary.

\subsection{Arcs, arc diagrams, and weighted arc-diagrams} \label{subsec:arcs}
An \emph{arc} in $S$ is a homotopy class of non-trivial \emph{embedded} paths (where we also mod out by reparametrization); 
we often write $\alpha = [a]$ for the arc $\alpha$ represented by the path $a$.
We say that two arcs are \emph{disjoint} if they have disjoint representatives.
In particular, an arc is disjoint from itself\footnote{We can rule out this case by saying that the arcs are \emph{distinct} and disjoint}. 
We let $\AAA(S)$ be the set of all non-trivial arcs on $S$.
If $\alpha, \beta \in \AAA(S)$, 
we let $\pair\alpha\beta$ be the minimal number of intersection points
 of representatives of $\alpha$ and $\beta$. 
In particular, 
$\pair \alpha\beta = 0$ if and only if $\alpha$ and $\beta$ are disjoint.  
We observe that an arc $\alpha$ has well-defined endpoints 
(in $\EE(S)$,
the set of endpoints of $S$), 
because the endpoints of a path are homotopy-invariant.

An {\it arc diagram} or \emph{multi-arc} on $S$ is a set  $\ba$ of pairwise disjoint arcs $\alpha_i$ (on $S$).
We observe that any arc diagram consists of at most $3|\chi(S)|$ distinct arcs.
A {\it weighted arc diagram (WAD) } (or weighted multi-arc) $X$ on $S$
is a map $X \from \AAA \to [0, \infty)$ such that any two arcs in the support of $X$ are disjoint.
(The support of $X$, denoted $\supp X$, is the arcs $\alpha$ for which $X(\alpha) > 0$). 
We will also think of $X$ as the formal weighted sum $\sum_\alpha X(\alpha) \alpha$ of arcs. 
We let $\WW(S)$ stand for the set of WAD's on $S$.
For $X=\sum w_i \alpha_i$ and $Y = \sum v_j\beta_j$,
we let 
$$\pair XY = \sum_{i, j} w_i v_j \pair{\alpha_i}{\beta_j}.$$
In particular, 
$\pair XY = 0$ if and only if no arc in $\supp X$ intersects an arc in $\supp Y$, 
and hence $\pair XX = 0$ for any WAD $X$.
We can think of arcs as arc-diagrams, and arc-diagrams as weighted arc-diagrams 
(with weight one on each arc in the diagram),
and thereby pair any two of the these objects. 

When $\pair XY = 0$ we can define $X+Y$ in the obvious way.
It is also useful to define subtraction for arc diagrams in terms of non-negative subtraction,
where
$
x \nnm y = \max(x - y, 0)$.
For any $X, Y \in \WW(S)$, we define $X-Y$ by $(X-Y)(\alpha)= X(\alpha) \nnm Y(\alpha)$.
Likewise, for $X \in \WW(S)$ and $c \in [0, \infty)$, 
we define $X-c$ by $(X - c)(\alpha) = X(\alpha) \nnm c$.

The set $\WW(S)$ is  {\it partially ordered}:
we say $X\leq Y$ if $X(\alpha)\leq Y(\alpha)$ for any $\alpha\in \supp X$.

We will make use of two norms on the space of WAD's:
$$
   \| W\|_\infty= \sup_{\alpha\in \AAA} W(\alpha); \quad \|W\|_1=\sum_{\alpha\in \AAA} W(\alpha).
$$

If $f\from U\ra V$ is a finite-sheeted  holomorphic covering between two Riemann surfaces
then there is a natural {\it pull-back} operation $f^*\from \WW(V)\ra \WW(U)$
acting on the WAD's.

\subsubsection{Valid WAD's} \label{subsubsec:valid}
A {\it proper lamination} $\FF$ on $S$ is a Borel set $Z\subset S$
explicitly realized as a union of disjoint proper paths called the {\it leaves} of $\FF$.
Any proper lamination\footnote{In what follows, all laminations under consideration are assumed to be proper.}
can be written $\FF = \bigcup_\alpha \FF(\alpha)$,
where $\FF(\alpha)$ comprises the leaves of $\FF$ that represent $\alpha$.
The arcs $\alpha\in \AAA$ for which $\FF(\alpha)$ is non-empty form an arc diagram.
For each such $\alpha$,
we can think of the sublamination $\FF(\alpha)$ as a path family on $S$,
and we let  $W_\FF(\alpha)$ be the \emph{extremal width}\footnote{Defined in Section \ref{A:def}}
$\WW(\FF(\alpha))$.
In this way we obtain the WAD
$W_\FF = \sum_\alpha W_\FF(\alpha)\alpha$
corresponding to $\FF$.

Note that if $f\from U\ra V$ is a finite-sheeted holomorphic covering between two Riemann surfaces
and $\FF$ is a proper lamination on $V$ then $f^*(\FF)$ is a proper lamination on $U$
and $W_{f^*(\FF)} = f^* (W_\FF)$.

Weighted arc diagrams that are $W_\FF$ for some proper lamination $\FF$ are called {\it valid}.

\subsection{The canonical foliation and weighted arc-diagram}
\label{subsec:canonical}
In this subsection,
and in the rest of the paper,
unless otherwise specified,
we will assume that the ideal boundary $\partial S$ is non-empty. 

Let us consider the universal covering $\pi\from \D\ra S$.
Let $\Gamma$ be the Fuchsian group of deck transformations of $\pi$,
and let $\Lambda\subset \T$ be its limit set. Since $\di S \neq \emptyset$,
$\Lambda$ is a Cantor set. 
We let $\hat S = \cl \D \sm \Lambda$,
so that  $\pi$ extends to $\pi\from \hat S \to S \cup \di S$ as the universal cover. 
We let $\di \hat S = \pi^{-1}(\di S) = \di \D \sm \Lambda$. 

Let us pick  two components, $I\not= J$, of $\di \hat S$.
The disk $\D$ with these two intervals as horizontal sides determines a quadrilateral $Q(I, J)$.
This quadrilateral can be conformally uniformized, $\phi\from  Q(I, J)\ra {\bf Q} (a) $,
by a standard quadrilateral ${\bf Q}(a) = [0,a]\times [0,1]$
in such a way that $I$ and $J$ correspond to the horizontal sides of ${\bf Q}(a)$.
The {\it vertical foliation} $ \FF(I, J)$ on $Q(I, J)$
is the $\phi$-pullback of the standard vertical foliation  on ${\bf Q}(a)$.
\begin{figure}[htbp]
\begin{center}
\makeabbrevs
\input{canonical2.pdf_t}
\caption{The canonical foliation}
\label{canonical}
\end{center}
\end{figure}

Assume now that $a>2$, and let us cut off from ${\bf Q}(a) $ two side squares, $[0,1]\times [0,1]$
and $[a-1, a]\times [0,1]$.  We call the left-over rectangle ${\bf Q}_\can(a)$,
and we let $Q_\can( I,  J) = \phi^{-1}({\bf Q}_\can(a))$.
The side quadrilaterals that we have cut off from $Q(I,J)$ are called its {\it buffers}.

Let $\FF_\can(I, J)$ be the restriction of $ \FF(I, J)$ to $Q_\can(I, J)$.
Obviously, for any deck transformation $\gamma \in \Gamma$, we have:

\begin{equation}\label{phi}
        \FF_\can ({\gamma (I), \gamma (J)}) = \gamma(\FF_\can (I, J)).
\end{equation}

\begin{lem}\label{disjointness}
  The rectangles $Q_\can (I, J)$ are  pairwise disjoint.
\end{lem}
\begin{figure}
\includegraphics[width=3in]{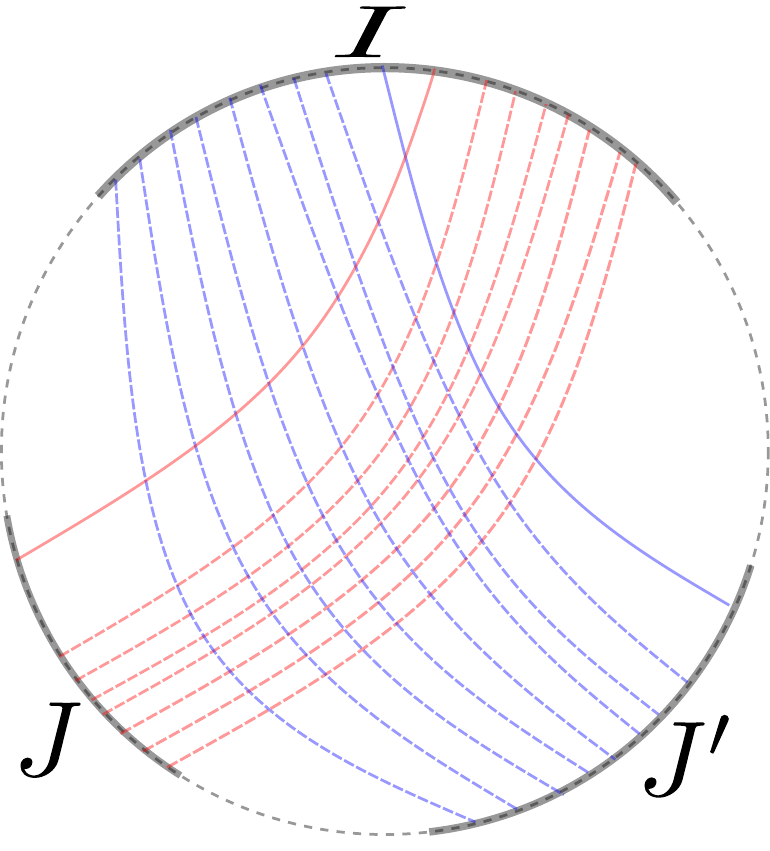}
\caption{Here a leaf of $\FF_\can(I, J)$ crosses a leaf of $\FF_\can(I, J')$. (These are the solid red and blue paths.)
This forces all paths in one buffer $B$ crossing all the paths of another buffer $B'$. (There are the dashed red and blue paths.)}
\label{fig:crossing}
\end{figure}
\begin{proof}
\def\Tt{T}
Let us consider two rectangles, $Q \equiv Q(I, J)$ and $Q' \equiv Q(I', J')$.
Then we can find one interval from each pair,
say $J$ and $J'$,
such that $J \neq J'$
(so $J \cap J' = \emptyset$),
and such that there is a component $\Tt$ of $\T - (J \cup J')$
such that $\Tt \cap (I \cup I') = \emptyset$.

We let $B$ be the buffer of $Q(I,J)$
 that has a horizontal side (which is a subset of $J$) that shares an endpoint with $\Tt$;
we define $B'$ likewise.
Then if any vertical leaf $\gamma$ of $\FF_\can(I, J)$
 crossed any vertical leaf $\gamma'$ of $\FF_\can(I', J')$,
then every vertical leaf of $B$
 would cross every vertical leaf of $B'$
 (see Figure \ref{fig:crossing}).
This would contradict Lemma \ref{no crossing}.
\end{proof}
This lemma allows us to define the {\it canonical lamination}  $\FF_\can(\hat S)$ as the union of
the laminations $\FF_\can(I, J)$
for all pairs of different components $I$ and $J$ of $\di\hat S$.
By (\ref{phi}), this lamination is $\Gamma$-invariant,
and hence it can be pushed forward to  $S$.
 In this way we obtain the {\it canonical lamination on $S$}:
$$
    \FF_\can(S) =  \pi_* ( \FF_\can(\hat S) ).
$$
The corresponding weighted arc diagram $\alpha\mapsto W_\can(S, \alpha) $, $\alpha\in \AAA(S)$, is  called
the {\it canonical WAD} on $S$.\footnote{We will use abbreviated notations $W_\can(S)$ or $W_\can(\alpha)$ whenever it does not lead to confusion.}
 By definition, it is valid.

We will now list several basic properties of the canonical WAD.
The rest of the theory will be based on these properties in an
essentially axiomatic way.

\subsection{Property A: Maximality}
For any arc $\alpha\in \AAA$,
we let $W_{\max}(S, \alpha) \in \R_+$
be the extremal width of the family of all proper paths $\gamma$  in $S$ representing $\alpha$.
We then have

\begin{lem}\label{max}
   For any valid arc diagram $W$ on $S$ and arc $\alpha\in \AAA(S)$, we have
$$
   W(\alpha) \leq W_{\max}(S, \alpha)\leq  W_\can(S, \alpha)  + 2.
$$
 \end{lem}

\begin{proof}
  The first inequality follows immediately from the definition of a valid arc diagram.
  
Let us prove the second inequality. 
  It is trivial for any arc $\alpha\in \AAA$ with $W_\max(\alpha)\leq 2$.
Let us consider some arc $\alpha\in \AAA$ with $W_\max(\alpha)>2$.
This arc connects two boundary components, $\si$ and $\om$, of $S$.

The path family $\GG(\alpha)$ representing $\alpha$
lifts to a path family $\hat \GG(\alpha)$
consisting of all the paths in $\D$ that connect
two appropriate arcs on $\T $,  $I$ and $J$,  covering  $\si$ and $\om$ respectively.
Viewing $ I$ and $ J$
as the horizontal sides of the rectangle $Q(I, J)$ based on $\bar \D$,
we obtain the desired estimate:
$$
       \WW(\GG(\alpha))\leq \WW(\hat\GG(\alpha))\leq \WW(Q(I,  J) ) = W_\can(S, \alpha)+2
$$
(where we have made use of Lemma \ref{increase} for the first estimate).
\end{proof}

\subsection{Property B: Natural behavior under coverings}

\begin{lem}\label{cov}
   If $f\from U\ra V$ is a finite-degree covering then $W_\can(U)= f^* W_\can(V)$.
\end{lem}

\begin{proof}
  Let $\pi_U\from \hat U\equiv \cl\D\sm \Lambda_U\ra U$ and $\pi_V\from \hat V\equiv \cl\D\sm \Lambda_V\ra V$
 be the universal  coverings of $U$ and $V$, with deck transformations groups $\Gamma_U$ and $\Gamma_V$ respectively.
Since $f$ has a finite degree,
the group $\Gamma_U$ has a finite index in $\Gamma_V$.
It follows that $\Lambda_U=\Lambda_V$, so that $\hat U=\hat V$.
Hence $\FF_\can(\hat U)=  \FF_\can (\hat V)\equiv \FF$.
Then
$$
 \FF_\can (U) =\FF/\Gamma_U = f^* (\FF/\Gamma_V)= f^* (\FF_\can(V)),
$$
and the conclusion follows.
\end{proof}

\subsection{Property C: Arcs with properly mapped endpoints}
Recall that $\EE(S)$ is the set of ends of the Riemann surface $S$.
Suppose that $A \in \EE(U)$ and $e\from U \to V$ is holomophic. 
We say that $e$ is \emph{proper at $A$} if $e(x_n) \to \infty$ in $V$ whenever $(x_n) \to A$ in $U$. 
In this case there is a unique $B \in \EE(V)$ such that $e(x_n) \to B$ whenever $x_n \to A$,
and we write $B = e(A)$. 
We call $A$ a \emph{proper end for $e$}. 

If $r$ is a path in $U$ that connects two proper ends of $e$ (where $e\from U \to V$ is as before),
then $e(r)$ is a path in $V$. 
Moreover homotopic paths are mapped to homotopic paths,
so $e(\rho)$ is defined for any arc $\rho$ on $U$ whose ends are proper for $e$. 
The only problem is that $\rho$ may not be embedded!
So let us now assume that $e$ is homotopic to an embedding;
then we can be sure that $e(\rho)$ is embedded.

We now state Property C for canonical WAD's:
\begin{lem}\label{e}
Suppose that $e\from U \to V$ is a holomorphic map between finite topology Riemann surfaces, 
and $e$ is homotopic to an embedding. 
Then for any arc $\alpha$ on $U$ whose endpoints are proper ends for $e$,
$$
     W_\can(U, \alpha)  \leq  W_\can(V)(e(\alpha)) .
$$
\end{lem}

\begin{proof}
The arc $\alpha$ connects two proper ends, $\si$ and $\om$,
 and it lifts to an arc in $\hat U$ connecting
some components $I$ and $J$ of $\di \hat U$.
Let $Q = Q(I, J)$  as defined in Section \ref{subsec:canonical}.
We observe that the lift of $e$ to the universal covers of $U$ and $V$ extends continuously to $I \cup J$;
we will denote this extension of the lift by $\hat e$. 
Let $I'=\hat e(I)$, $ J'=\hat e(J)$,
and let $Q' = Q(I', J')$.
Then $\hat e\from Q\ra Q'$ maps the horizontal sides of $Q$ to the corresponding horizontal sides of $Q'$.
By Corollary \ref{Q}, $\WW(Q)\leq \WW(Q')$. But $\WW(Q)= W_\can (U, \alpha)+2$, $\WW(Q')= W_\can (V, e_*(\alpha))+2$,
and the desired conclusion follows.
\end{proof}

We can also state the following corollary,
which we will also call Property C.
With $e\from U\to V$ as before,
let $\EE_e(U)$ be the proper ends for $e$,
and let $\EE_e(V) = e(\EE_e(U))$. 
Let $\AA'(U)$ be the arcs with endpoints in $\EE_e(U)$,
let $\Wcan'(U) = \Wcan(U)|_{\AA'(U)}$,
and define $\AA'(V)$ and $\Wcan'(V)$ analogously.
Then the map $e\from \AA'(U) \to \AA'(V)$ induces a pullback $e^*$ 
from weighted arc-diagrams supported on $\AA'(V)$ to WAD's supported on $\AA'(U)$. 
From Lemma \ref{e} we can immediately conclude that
\begin{equation}
\Wcan'(U) \le e^*\Wcan'(V). 
\end{equation}

\subsection{Property D: Domination}

Let us now introduce an important relation between WAD's.
We consider two Riemann surfaces, $U\subset V$.
Given a path $\gamma$ on $V$,
the restriction $\gamma\cap U$
has only finitely many non-trivial components,
$(\gamma_i)_{i=1}^n$.
They represent a sequence of arcs on $U$,
$I(\gamma) \equiv (\alpha_i\equiv[\gamma_i])_{i=1}^n$,
called the {\it itinerary} of $\gamma$.

We say that a sequence of arcs $(\alpha_i)$ on $U$ {\it arrows} an arc $\beta$ on $V$
if there exists a path $\gamma$ representing $\beta$ such that $I(\gamma)=(\alpha_i)$.
We will use notation $(\alpha_i) \arrow \beta$ for the arrow relation.

\ssk{\it Remark.}
   Note that the endpoint of $\gamma_i$ is connected to the beginning of $\gamma_{i+1}$
by a path that goes through some component $K$ of $V\sm U$. In this case,
the end of $U$ corresponding to this component is {\it not properly embedded}
into $V$. This remark is useful as it reduces a number of possibilities of
how the arc $\beta$ can be composed by arcs $\alpha_j$.
\ssk

Let us now consider two WAD's,  $X\in \WW(U)$ and $Y\in \WW(V)$.
We say that $X$ {\it dominates} $Y$,
written
$$X\multimap Y,$$
if we can write
$$
X\geq \sum_i \sum_j w_{ij}\alpha_{ij}, \quad Y = \sum_i v_i \beta_i
$$
where, for each $i$,
$$
(\alpha_{ij})_j \arrow \beta_i
$$
and
$$
\bigoplus_j w_{ij} \ge v_i.
$$
We develop the basics of the theory of domination in Appendix B. 

We have defined domination so as to have the following statement:
\begin{lem}\label{XY}
  Given a valid WAD $Y$ on $V$, there exists a valid WAD $X$ on $U$ such that $X\multimap Y$.
\end{lem}

\begin{proof}
Since $Y$ is valid,   $Y=W_\FF$ for some  lamination $\FF$ on $V$.
For $\beta\in \supp W_\FF$,
let $\FF(\beta)$ be the sublamination of $\FF$ assembled by the leaves $\gamma$
representing the arc $\beta$.

Let us consider the restriction of $\FF$ to $U$, that is, let $\HH = \FF\cap U$ and $X=W_\HH$.
To any leaf $\gamma$ of $\FF$, let us associate its itinerary $I(\gamma)=(\alpha_j(\gamma))$
on $U$.
Let $\II(\beta)$ stand for the set of all non-trivial itineraries $\ba = I(\gamma)$
corresponding to all possible leaves $\gamma$ of $\FF(\beta)$.
By definition, $ \ba \arrow \beta$ for any $ \ba\in \II(\beta)$.
Let $\FF(\beta, \ba)$ stand for the sublamination of $\FF(\beta)$ assembled by the leaves $\gamma$
with itinerary $\ba$, i.e., for which $I(\gamma) = \ba$.

For $\ba=(\alpha_j)\in \II(\beta)$,
let $v(\beta,\ba) =  W(\FF(\beta,\ba))$, and
let $w_j (\beta, \ba)$ be the width of the lamination assembled by the segments
of $\FF(\beta, \ba) \cap U$ corresponding to $\alpha_j$.
By the Series Law,
$$
   \bigoplus_j w_j (\beta, \ba)  \geq  v(\beta, \ba).
$$
  Moreover,
$$
     X =  \sum_\beta \sum_{\ba\in \II(\beta)} \sum_j  w_j(\beta, \ba) \alpha_j, \quad
     Y = \sum_{\beta} \sum_{\ba\in \II(\beta)} v (\beta, \ba) \beta.
$$
This would mean that  $X\multimap Y$ if we knew that $\II(\beta)$ were finite. 
For any $\gamma \in \FF$, 
we let $| I(\gamma) |$ denote the length of $I(\gamma)$, i.e. the length of the sequence $(\alpha_j(\gamma))$. 
For $n \in \N$, we let $\FF_n$ be the union on $\FF(\beta, \ba)$ over all $\beta$, and all $\ba$ with $|\ba| \le n$;
we let $Y_n = W_{\FF_n}$. 
Then $Y_n \to Y$ and $X \doms Y_n$ for all $n$, so $X \doms Y$ by Lemma \ref{lem:domination-closed}.
\end{proof}

\begin{rem}
It is not so very hard to show that $|\II(\gamma)|$ varies upper semi-continuously with $\gamma$ (in each $\FF(\beta)$).
Therefore $|\II(\gamma)|$ is bounded, and $Y = Y_n$ for some $n$. 
\end{rem}

We can now prove Property D of canonical WAD's:

\begin{lem}\label{B}
  Let $U\subset V$. Then there exists a WAD $B\in W(U)$ with $\|B\|_\infty\leq 2$ such that
$$
  W_\can(U) + B \multimap W_\can(V).
$$
\end{lem}

\begin{proof}
  Since $W_\can(V)$ is valid,
 Lemma \ref{XY} gives us a WAD $X$ on $U$ such that $X\multimap W_\can(V)$.
 By the Maximality Property, $X\leq W_\can(U)+2$. Hence there exists a WAD $B$ with
 $\supp B\subset \supp X$, $\|B\|_\infty\leq 2$,  and $X\leq W_\can(U)+B$. The conclusion follows.
\end{proof}

This implies the following, which we will also call Property D. 

\begin{cor}
   Let $U\subset V$. Then
$$
  W_\can(U)  \multimap W_\can(V)- 6 |\chi(U)|.
$$
\end{cor}
\begin{proof}
This follows immediately from Lemma \ref{B}, Lemma \ref{lem:subtract-from-right}, and the estimate $\norm{B}_1 \le 3|\chi(U)| \norm{B}_\infty$.
\end{proof}
\subsection{Property E: Intersection number and hyperbolic length}
Let $S$ be a compact Riemann surface with boundary. 
We have defined $\Wcan(S)$, 
and we have also defined the the weighted arc-diagram $M^S$ in Appendix C.
We can relate $M^S$ to $\Wcan^S$ via the following:

\begin{lem}\label{hyp dist vs mod}
  Let $Q$ be a quadrilateral with horizontal sides $I_1$ and $I_2$,
endowed with the hyperbolic metric.
Let $\Gamma_1$ and $\Gamma_2$ be the hyperbolic geodesics in $\Pi$
with the same endpoints as $I_1$ and $I_2$ respectively.
Then
$$
    \WW(Q) = -\frac{2}{\pi} \log \dist (\Gamma_1, \Gamma_2) + O(1).
$$
\end{lem}

\begin{proof}
We can map $Q$ conformally to the infinite strip
$$
\Pi \equiv \{ z : 0 < \Im z < \pi \}
$$
such that $\Gamma_1$ and $\Gamma_2$ map to vertical transverse segments
 with real part 0 and $d$ respectively,
where $d \equiv \dist(\Gamma_1, \Gamma_2)$.
Then the vertical sides of $Q$ map to
 $J_1 \equiv [0,d] \times \{0\}$
 and
 $J_2 \equiv [0, d] \times \{\pi\}$.
Then $\WW(Q)$ is equal to $\LL(\Theta)$,
where $\Theta$ is the family of paths in $\Pi$ connecting $J_1$ and $J_2$.
By two applications of the reflection principle (see \cite{Ahlfors:quasiconformal}),
we find that $\LL(\Theta) = 4 \LL(\Theta')$,
where $\Theta'$ is the family of paths connecting $J_1$ to the boundary of
$\Pi' \equiv \{ z : |\Im z| < \pi/2 \}$.
By the round annulus theorem\cite{towers}
$$
\LL(\Theta') = \frac{1}{2\pi} \log \frac {\pi/2}d + O(1) = - \frac{1}{2\pi}\log d + O(1)
$$
when $d$ is bounded above.
The theorem follows.
\end{proof}

We then have the following corollary:

\begin{lem} \label{two WAD's}
For any arc $\alpha \in \AA(S)$,
$$
|W_{\mathrm can}^S(\alpha) -\frac{2}{\pi} M^S(\alpha)| < C_0.
$$
\end{lem}

\begin{proof}
Given $\alpha$,
let us consider a lift $\tilde \alpha$ to the universal cover.
It connects two arcs $I_1$, $I_2$ on the circle
that cover the boundary curves of $S$ that the endpoints of $\alpha$ lie on.
We can lift $h\from \mathbf S \to S$ to $\tilde h\from \core R \to \D$,
where $\core R$ is the convex hull of the limit set of the deck transformation group.
Then letting $\Gamma_1$, $\Gamma_2$ connect the endpoints
 of $I_1$, $I_2$ respectively,
we find that the transverse geodesic arc $\bfalpha$
 lifts to the common perpendicular segment to $\Gamma_1$ and $\Gamma_2$.
So $\dist(\Gamma_1, \Gamma_2) = L(\bfalpha)$,
and the Lemma follows from Lemma \ref{hyp dist vs mod}.
\end{proof}
From Theorem \ref{peripheral lengths} and Lemma \ref{two WAD's} we can immediately conclude:
\begin{corollary} \label{gamma-W}
Let $S$ be a compact Riemann surface with boundary,
and endow $\Int S$ with its Poincar\'e metric.
Let $\gamma$ be a peripheral closed geodesic in $\Int S$.
Then
$$
L(\gamma) = \pi \pairing{ \Wcan(S), \gamma } + O(1; \chi(S)).
$$
\end{corollary}
With $S$ as in Corollary \ref{gamma-W},
we let $\Gamma_S$ be the formal sum of the peripheral closed geodesics in $\Int S$.
Then
$$
\pairing{ \Wcan(S), \Gamma_S } = 2 \norm{\Wcan(S)}_1,
$$
because every arc in $S$ intersects exactly two peripheral geodesics exactly once (or intersects exactly one exactly twice). 
Therefore,
by Corollary \ref{gamma-W},
\begin{equation} \label{eq:vh-gamma}
\norm{\wcan(S)}_1 = \frac\pi 2 |\Gamma| + O(1; p). 
\end{equation}
We can use this equation to compare the total weight of $\wcan(S)$ between two surfaces:
\begin{theorem} \label{thm:total-increase}
Suppose $S$ and $S'$ are finite topology Riemman surfaces,
and there is a holomorphic map $i \from S' \to S$ that is homotopic to a homeomorphism.
Then
$$
\norm{\wcan(S)}_1 \le \norm{\wcan(S')}_1 + O(1; \chi(S)).
$$
\end{theorem}
\begin{proof}
Because $i$ is homotopic to a homeomorphism,
it bijectively maps the free homotopy classes of the peripheral geodesics of $S'$ to the same for $S$.
Then, 
by the Schwarz Lemma,
because $i$ is holomorphic,
we have $|\Gamma_{S'}| \ge |\Gamma_S|$.
The Theorem then follows from the application of \eqref{eq:vh-gamma} to $\Gamma_{S'}$ and $\Gamma_S$. 
\end{proof}

\section{Theorem \protect{\ref{thm:length-comparison}} in terms of the canonical weighted arc-diagram} 
\label{sec:reduce-to-wcan}
In this brief section we will prove Theorem \ref{thm:length-comparison}, given a simply stated theorem about $\Wcan(V \sm \KK_p)$ in that context. 
We then outline the proof of this latter theorem, which is then proven in Sections \ref{Hubbard trees}--\ref{sec:hor-vert}.
 
 \subsection{Horizontal and vertical arc diagrams} \label{subsec:hv}
 
In further applications, the Riemann surface $S$ we worked with in Section \ref{sec:canonical} will be $V\sm \KK$,
where $V$ is a topological disk, and $\KK$ is the union of $p$ disjoint FJ-sets
(the cycle of small Julia sets of some quadratic-like map).

Under these circumstances,
a proper path (and the corresponding arc) in $V \sm \KK$ is called {\it horizontal} if it connects two (not necessarily distinct) components of $\KK$,
and is called {\it vertical} if it connects a component of $\KK$ to $\di V$.
Given an arc diagram $W$ on $V\sm \KK$, let $W^\hor$ and $W^\ver$  stand for its horizontal and vertical parts respectively,
and let $W^{\ver+\hor}$ stand for their sum.
In particular, we can consider  $\Wcan^\hor$, $\Wcan^\ver$ and $\Wcan^{\ver+\hor}$ for $V \sm \KK$.

We can now state a theorem which will be the main result of Section \ref{sec:hor-vert}:
\begin{theorem}\label{def vertical teaser}
Let $f\from U'\to U$ be primitively renormalizable quadratic-like-map with period $p$,
and let $\KK$ be its cycle of small Julia sets.
Then there exists $M=M(p)$ such that:
If 
\begin{equation} \label{eq:teaser:hyp}
\| W_\can^{\ver+\hor}(U\sm \KK)\|_1 > M
\end{equation}
then
\begin{equation} \label{eq:teaser:conc}
 \|W^\ver_\can(U\sm\KK)\|_1 \geq C^{-1} \| W^\hor_\can(U\sm \KK)\|_1,
\end{equation}
with an absolute constant $C>0$.
\end{theorem}

Given Theorem \ref{def vertical teaser}, we can prove Theorem \ref{thm:length-comparison}.
\begin{proof}[Proof of Theorem \ref{thm:length-comparison}]
In addition to the notation of Section \ref{intro:defs},
we let $\Gamma$ be peripheral closed geodesic in $U\sm \KK_p$ homotopic to $\di U$. 
We let $W_\can \equiv W_\can(U\sm \KK_p)$.

The geodesic $\Gamma$ intersects each vertical arc once and does not intersect horizontal arcs.
By  Corollary \ref{gamma-W},
\begin{equation}\label{Gamma}
  |\Gamma| \ge \pi \| W_\can^\ver\|_1 +O(1).
\end{equation}

 Let $W_\can|_j \subset W_\can^{v+h}$ be the part of $W_\can$
supported on the arcs landing at $K_p(j)$. The geodesic $\gamma_p(j)$ does not intersect
$W_\can - W_\can|_j$ and intersects each arc of $\supp W_\can|_j$ once or twice.
By  Corollary \ref{gamma-W},
\begin{equation}\label{gamma_j}
   \pi\| W_\can|_j\|_1 + O(1) \le |\gamma_p(j)| \le 2\pi\| W_\can|_j\|_1 + O(1).
\end{equation}
But by the Schwarz Lemma, the geodesics $\gamma_p(j)$ have comparable lengths
(see \cite[Theorem 9.3]{McM}):
\begin{equation} \label{eq:compare}
   \frac{1}{2} |\gamma_p| \leq |\gamma_p(j)| \leq |\gamma_p|.
\end{equation}
Putting this together with (\ref{gamma_j}), we see that all the $\| W_\can|\, j \|_1 $
are also comparable, provided $|\gamma_p|$ is sufficiently big
(bigger than some absolute $L_0$). Hence
$$
   \| W_\can|_0\|_1  \asymp \frac{1}{p} \| W_\can^{\ver+\hor}\|_1,
$$
and together with (\ref{gamma_j}), we obtain:
$$
     |\gamma_p| \asymp \frac{1}{p} \| W_\can^{\ver+\hor}\|_1.
$$

Now assume $ \| W_\can^{\ver+\hor}\|_1 > M(\bar p)$. 
Then Theorem \ref{def vertical teaser} implies that $\| W_\can^{v+h}\|_1$ is comparable with $\| W_\can^\ver\|_1 $,
and we can conclude
$$
   | \gamma_p | \leq \frac{ C }{p} \| W_\can^\ver\|_1,
$$
which, 
combined with \eqref{Gamma},
implies
$$
   |\gamma_p| \leq \frac{C}{p} |\Gamma|.
$$
We now need only observe that by the Schwarz Lemma, $|\Gamma| \leq  |\gamma|$ (since $K\supset \KK_p$).
\end{proof}

\subsection{Outline for the proof of Theorem \ref{def vertical teaser}} \label{subsec:hv-proof-outline}
In this informal subsection we will outline the proof of Theorem \ref{def vertical teaser} that will be presented in Sections \ref{sec:entropy} and \ref{sec:hor-vert}, using the background in Section \ref{Hubbard trees}.
For the sake of clarity we have slightly simplified the argument and altered some of the constants.

\subsubsection{Canonical WAD's in restricted domains}
In Section \ref{sec:hor-vert},
we will consider the nest of the domains $U^m := f^{-m}(U)$,
and their associated canonical WAD's $\Wcan(U^m \sm \KK_p)$.
Let $X^m$ stand for the horizontal parts of these WAD's
(with appropriate constants subtracted).
We show that 
\begin{enumerate}
\item
$f^* X^m\multimap X^{m+1}$ and 
\item  $X^{m+1} \le X^m$. 
\end{enumerate}

\subsubsection{Loss of horizontal weight through alignment and entropy}
In Section \ref{sec:entropy},
after a discussion in Section \ref{Hubbard trees} of some dynamical properties of the Hubbard tree, 
we show the purely combinatorial theorem\footnote{While the motivation for this combinatorial theorem relies on the construction of the $X^n$ in Section \ref{sec:hor-vert}, the proof does not.}, 
that Properties 1 and 2 of $X^m$ above imply that 
\begin{equation} \label{eq:prelim:reduction}
\norm{X^{8p}}_1 \le  \frac 14\norm{X^0}_1.
\end{equation}
This proceeds in two stages: we first show that $X^{3p}$ is aligned with the Hubbard tree 
of the superstable model for $f$. 
We then observe that dynamics on the Hubbard tree is expanding (has positive entropy);
this, 
coupled with the alignment of $X^{3p}$ with the Hubbard tree,
allows us to show \eqref{eq:prelim:reduction}.
This is the main place where we use that the renormalizations are
of primitive type.

Returning to Section \ref{sec:hor-vert},
we see that the application of the combinatorial statement \eqref{eq:prelim:reduction} to our geometric situation implies,
assuming \eqref{eq:teaser:hyp},
that
 $$
     \|W_\can^\hor(U^{8p}\sm \KK_p)\|_1 \leq \frac{1}{3} \| W^\hor_\can(U^0 \sm \KK_p) \|_1.
$$

\subsubsection{Push-forward via the Covering Lemma}
   The horizontal width released under restrictions is turned into the vertical  width (Section \ref{subsec:vert-from-hor-loss}),
which implies
$$
 \norm{W_\can(U^{8p} \sm\KK_p)}_1\geq  \frac{1}{2}\norm{W_\can^\hor(U^{8p}\sm \KK_p)}_1.
$$
We can then apply the Covering Lemma\cite{covering lemma} in Section \ref{covering-application} to push forward $W_\can(U^{8p} \sm\KK_p)$  by  $f^{8p}$
and obtain \eqref{eq:teaser:conc}.

\section{Some properties of Hubbard trees}\label{Hubbard trees}
In this section we will introduce the superstable model $F$ for a renormalizable map $f$ of period $p$.
We will then relate the Hubbard tree for $F$ to arcs of the kind we studied in Section \ref{sec:canonical},
and show that it has certain simple dynamical properties. We will then use these properties in Section \ref{sec:entropy} to complete the proof of Theorem \ref{loss-F}.

\subsection{The superstable model} \label{subsec:model}
\newcommand{\cdd}{{\DD}}

Suppose that $f\from \C\to \C$ is a given by $f(z) = z^2 + c$,
and $f$ is primitively renormalizable of period $p$ and $\KK\equiv \KK_p$ is the union of small Julia sets of $f$.
Then $c$ belongs to a \emph{small Mandelbrot set} $M_{c_0}$,
where $F(z) = z^2 + c_0$ has an superattracting periodic cycle of period $p$.
We call $F$ the {\it superstable model}  for $f$.
Let $\OO=\{F^k (0)\}_{k=0}^{p-1}$ be the superattracting cycle of $F$,
and $\DD$ be the closure of the immediate attracting basin for the superattracting cycle for $F$;
it is also the union of small Julia sets of the primitively $p$-renormalizable map $F$.
Then we can write $\DD = \bigcup_{j=0}^{p-1} D_j$,
where $D_j$ is the component of $\DD$ containing $F^j(0)$.
We observe that the $D_j$'s are disjoint closed Jordan domains.

The proof of the following lemma is left to the reader;
it can be proven using the similarity of the Yoccoz puzzles of $f$ and $F$,
or using Thurston's theorem \cite{DH-Thurston},
and the collapsing of the $\KK$, as in \cite[Appendix B]{McM}.
\begin{lemma}
There is a homeomorphism $h\from \C\sm \cdd \to \C \sm \KK$ that  lifts to 
$\tilde h \from \C \sm F^{-1} \cdd \to \C \sm f^{-1} \KK$ such that
$h \sim \tilde h$ on $\C \sm F^{-1} \cdd$. 
Moreover, 
$F$ and $h$ (up to isotopy) are characterized by this property.
\end{lemma}
We call $h$ a \emph{Thurston equivalence} between $F$ and $f$, because it is essentially the same as the object of the same name described in \cite{DH-Thurston}.

Similarly,
if $f\from U' \to U$ is a $p$-renormalizable quadratic-like map (with union of small Julia sets $\KK$),
then we can find a critically periodic $F$ as above and a Thurston equivalence $h\from \C \sm {\DD} \to U \sm \KK$
(lifting to $\tilde h \from \C \sm F^{-1}  \DD \to U' \sm f^{-1} \KK$ such that
$h \sim \tilde h$ on $\C \sm F^{-1} \cdd$).

If $\alpha$ is an arc on $U \sm \KK$,
we can \emph{transfer} $\alpha$  to the arc $h^{-1}(\alpha)$ on $\C \sm \cdd$. 
We can likewise transfer arc-diagrams and weighted arc-diagrams;
we observe that all the combinatorial relations studied in Section \ref{sec:canonical}
are preserved under transfer.

\subsection{The Hubbard tree}
In the remainder of this section we will work exclusively with the superstable model $F$,
with the ultimate aim of proving certain combinatorial properties for the dynamics of arcs on $\C \sm \DD$.

Given a set $X\subset K(F)$, the {\it path hull} $\Hull(X)$ of $X$ is defined as the smallest
path connected closed subset of $K(F)$ containing $X$ such that $\Hull(X) \sm X$ intersects any component
of $\inter K(F)$ in the union of internal rays.  The {\it Hubbard tree} $H\equiv H_F$
is defined as the path hull of the basin $\DD$.
It is invariant under $F$; in fact, $F(H)=H$.
We observe that $H$ is simply connected, and therefore 
the closure of every component of $H \sm \cdd$ is an embedded and simply connected 1-complex in $\C \sm \DD$
that intersects each $D_j$ in at most one point.

If $D$ is a component of $\DD$,
the \emph{valence} of $D$ (in $H$) is the number of components of $H \sm  \DD$ whose closures intersect $D$. 
The \emph{cut number} of $D$ is the number of components of $H \sm  D$. 
We observe that the valence of each $D$ is equal to its cut number.

We close with one lemma:
\begin{lem}\label{valence}
The valence of $D_1$ is 1 and the valence of $D_p = D_0$ is 2. 
The valence of $D_i$ is non-decreasing in $i$ for $i = 1, \ldots, p$. 
Therefore the valence of any disk $D_j$ of the Hubbard tree is at most~2.
\end{lem}
\begin{proof}
Because $F$ is locally injective at $D_i$, for $0 < i < p$,
the valence of $D_i$ is at most that of $D_{i+1}$. 
Similarly the valence of $D_0$ is at most twice that of $D_1$.
But some $D_i$ has valence 1 because the Hubbard tree is a tree. 
It follows that $D_1$ has valence 1, and $D_0$ has valence at most 2.
But some $D_i$ has valence 2, because the Hubbard tree is connected and non-trivial. 
Therefore $D_p = D_0$ has valence 2.  
\end{proof}

\subsection{Disked trees, aligned arcs, and trees of complete graphs}
We will now properly define a disked tree, and then define the graph of aligned arcs in a disked tree, and see that it is a tree of complete graphs. 

\begin{defn}
A \emph{disked tree} $H$ is a connected, locally connected, and simply connected compact subset of $\C$, 
such that there is a set $\{ D_i \}$ of disjoint closed Jordan domains (with union  $\DD = \bigcup_i D_i$), 
such that $\DD \subset H$ and $H \setminus \DD$ has finitely many components $H_k$,
and the closure of each $H_k$ is an embedded finite 1-complex. 
\end{defn}
We call the $D_i$'s the \emph{disks} of $H$.
We observe that the Hubbard tree $H_F$ of our critically periodic $F$ is a disked tree, with disks $\{D_j\}_{j=0}^{p-1}$.

\begin{figure}
\includegraphics{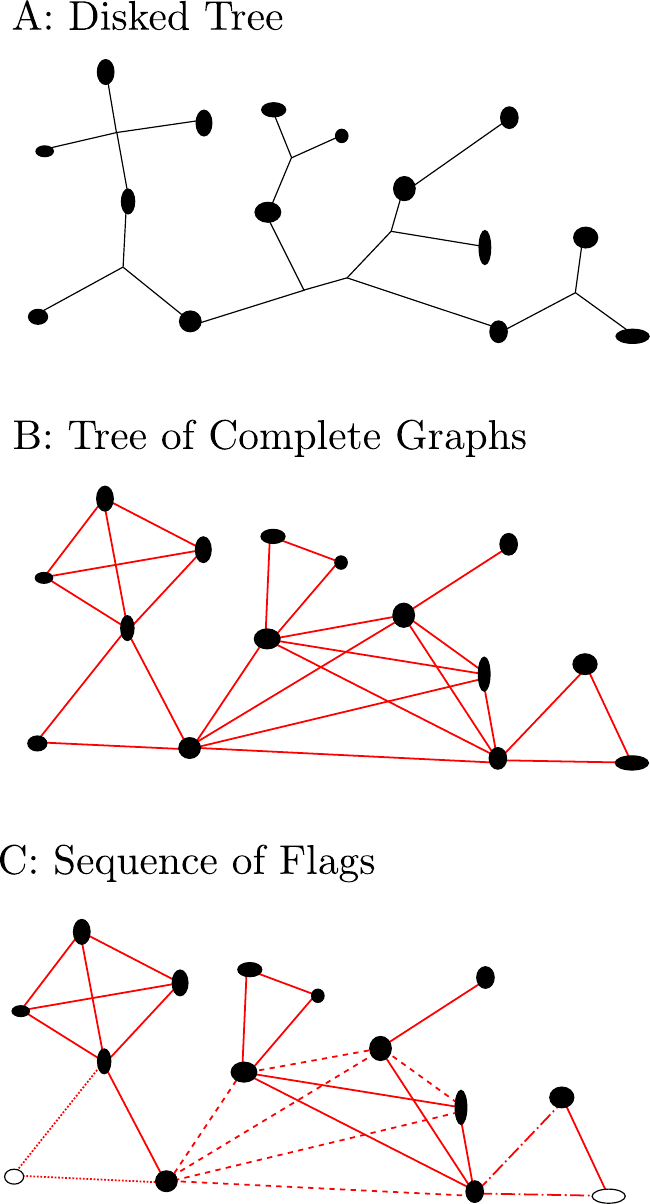}
\caption{A: We are given a disked tree. B: We can form the associated tree of complete graphs. C: Any path (of red arcs) between the two hollow disks must include at least one each of the dotted, dashed, and dashed-and-dotted arcs. This illustrates Lemma \ref{lem:into-flags}.}
\label{fig:variants}
\end{figure}

We say that a (proper) injective path in $\C \sm \DD$
is {\it aligned} with a disked tree $H$
if it lies in one of the 1-complexes $H_k$.
We say an arc in $\C \sm \DD$ is aligned with $H$ if it is represented by a path that is aligned with $H$.
Let $\BH$ stand for the family of the aligned paths/arcs.
(Since there is a natural one-to-one correspondence between the aligned arcs and
the aligned paths, we will not distinguish notationally these families).
We can think of $\BH$ as a graph where the vertices are the $D_i$ and the edges are these aligned arcs. 
See Figure \ref{fig:variants} for an example.

For each component $H_k$ of $H$, we let $\partial H_k$ be the set of $D_i$'s on which an endpoint of $H_k$ lies.
Because $H_k$ is path connected, 
every pair of $D_i$'s in $\partial H_k$ can be connected by a path lying in $H_k$, 
and hence in $\BH$. 
Because $H_k$ is simply connected, 
there is a unique such path,
so there is a unique arc in $\BH$ that connects this pair.
In this way we see that $\BH$ canonically contains a complete graph on $\partial H_k$;
we also see that it is the union (over all $k$) of these complete graphs.

The simple connectivity of $H$ also implies that for any disk $D$,
the number of complete graphs (in $\BH$) that $D$ belongs to is equal to the number of components of $\BH - D$, 
with a natural one-to-one correspondence between them.
We will call a union of edge-disjoint complete graphs with this property a \emph{tree of complete graphs}, in analogy to a tree, 
where for each vertex $v$ of the tree, there is exactly one component of the tree minus $v$ for each edge with $v$ as an endpoint.

Just as we can build every tree by starting with graphs consisting of a single edge connecting two vertices and then repeatedly joining disjoint graphs at a vertex (of each),
we can build every tree of disjoint graphs by starting with a collection of complete graphs and repeatedly joining disjoint graphs at a vertex. 

Suppose $\BG$ is a tree of complete graphs,
$v$ is a vertex of $\BG$, and $G \subset \BG$ is a complete graph in $G$ that has $v$ as a vertex.
Then the \emph{flag} $\Flag(v, G)$ is the set of edges of $G$ that have $v$ as an endpoint.
We observe that if $e$ is an edge of $\BG$ that has $v$ as an endpoint, 
there is a unique complete graph $G(e)$ of $\BG$ such that $e$ is an edge of $G$.
So we let $\Flag(v, e)$ be $\Flag(v, G(e))$. 

We then observe the following, which is left as a simple exercise for the reader.
See Figure \ref{fig:variants} for an illustration.
\begin{lemma}\label{lem:into-flags}
Suppose $\BG$ is a tree of complete graphs,
and $v$ and $w$ are distinct vertices of $\BG$.
Then there is a unique vertex-injective path in $\BG$ connecting $v$ and $w$.
Denoting the vertices of this path by $v_0 = v, v_1, \ldots v_n = w$,
and the edges by $e_0, e_1, \ldots, e_{n-1}$ (so each $e_i$ connects $v_i$ and $v_{i+1}$),
we find that the flags $\Flag(v_i, e_i)$ are disjoint (for $i = 0,\ldots, n-1$),
and every path connecting $v$ and $w$ has an edge in each of these flags. 
\end{lemma}

\subsection{Lifts of horizontal and vertical arcs}
For the rest of this section,
we return to the setting of a superstable quadratic polynomial $F$,
along with its disked tree $H$ and associated tree of complete graphs $\BH$. 
As we defined in Section \ref{subsec:hv},
a horizontal path or arc for $\C \sm \DD$ is one that connects two components of $\DD$,
and a vertical path or arc connects a component of $\DD$ to $\infty$. 
If a vertical path connects $D \in \DD$ to $\infty$, we will say that the path \emph{lands} at $D$. 

Suppose $\gamma_0$ and $\gamma_1$ are (by definition proper) paths in $\C \sm \DD$ and $F(\gamma_1) = \gamma_0$.
We say then that $\gamma_1$ is a \emph{proper lift} of $\gamma_0$.
If $\alpha_0 = [\gamma_0]$ and $\alpha_1 = [\gamma_1]$ are the associated arcs,
we say that $\alpha_1$ is a lift (by $F$) of $\alpha_0$. 
(Another point of view is that $\gamma_0$ has two lifts that are paths in $\C \sm F^{-1}\DD$,
but we keep only the lifts that have endpoints in  $\DD$.)
We observe that if $\alpha_0$ is horizontal, it can have only one lift; it may also have no lifts. 
If $\alpha_0$ is vertical, 
then it has exactly one lift, 
unless it lands at $D_1$, 
in which case it has two.

\subsection{An expanding metric}
Recall that $\OO = \{ F^i(0)\}_{i = 0}^{p-1}$ is the critical cycle for $F$.
We let $\rho_F$ be the Poincar\'e metric on $\C \sm \OO$.
We observe that $F\from \C \sm F^{-1}(\OO) \to \C \sm \OO$ is expanding for this metric,
and uniformly expanding on compact sets.
In particular,
for any bounded set $E$,
we see that $F$ is uniformly expanding on $E \sm F^{-1} \DD$. 
If $p$ is a rectifiable path in $\C \sm \OO$,
we let $|p|$ be the length of $p$ in this metric. 

We can then conclude:
\begin{lemma} \label{lem:contract-lifts}
Suppose that $p$ is a rectifiable path in $\C \sm \DD$ with (at least) one endpoint in $\cl\DD$,
and $q$ is a lift of $p$, in the sense that $F(q) = p$.
Then 
\begin{equation}
|q| \le \lambda |p|,
\end{equation}
where $\lambda < 1$ depends only on an upper bound for $|p|$. 
\end{lemma}

\subsection{No periodic horizontal arcs}
An arc $\alpha$ is called {\it periodic} if there is a sequence $\alpha_0, \ldots, \alpha_l$ such that $\alpha_l = \alpha_0$,
and each $\alpha_{i+1}$ is a lift of $\alpha_i$. 

We then have the following. 
\begin{lem}[see \cite{P}, Theorem 5.8]\label{no inv hor arcs}
 Periodic horizontal arcs do not exist.
\end{lem}

\begin{proof}
For every horizontal arc $\alpha$,
we let $|\alpha|$ be the infimum of lengths of paths representing $\alpha$,
(with respect to the Poincar\'e metric on $\C\sm \OO$).
It follows from Lemma \ref{lem:contract-lifts}
that if $\alpha'$ is a lift of $\alpha$,
then $|\alpha'| < |\alpha|$. 
Then given a cycle $\alpha_0, \ldots, \alpha_l = \alpha_0$ of lifts, 
we have an obvious contradiction. 
\end{proof}

\subsection{Trees $H^l$ and the associated objects.}\label{H-l}
For any  $l\in \N$, we can define the following objects:
\begin{itemize}
\item
  $\DD^l= F^{-l}(\DD)$;
\item
$H^l=F^{-l}(H)$, a disked tree whose disks are the components of $\DD^l$;
\item
$\BH^l$, the tree of complete graphs formed by the paths (or arcs) aligned with $H^l$. 
\end{itemize}
We observe that $H^{l+k} = F^{k*}H^l$ and $\BH^{l+k} = F^{k*}\BH^l$ for all $k, l \in \N$.
We also observe that every path $\gamma$ in $\BH^l$ is a uniquely determined concatenation of paths $(\gamma_i)$ in $\BH^{k+l}$,
in the sense that the intersection of $\gamma$ with $\C - \DD^{l +k}$ is the union of the (mutually disjoint) $\gamma_i$'s. 

Finally,
we observe that every complete graph in $\BH^l$ maps bijectively (by $F^l$) to a complete graph in $\BH$, and hence every flag $\Flag(v, e)$ in $\BH^l$ maps bijectively to $\Flag(F^l(v), F^l(e))$, a flag in $\BH$. 

\subsection{Expansion property}
  Let us consider a $|\BH| \times |\BH| $ matrix $M=M(F)$ defined as follows:
$M_{\gamma\de}=n$ if $\gamma$ contains $n$ disjoint segments $\gamma_i\in F^{-1}\BH$
such that  $F\gamma_i=\de$. One can easily check that, 
applying the same construction to the map $F^k$,
we obtain the matrix $M^k$.

\begin{lem}\label{expansion}
   For any $\gamma\in \BH$, we have $\displaystyle{\sum_\de M^p_{\gamma\de}\geq 2}$.
\end{lem}

\begin{proof}
If not,
then $\gamma$ (thought of as a path) would intersect no disk in $\DD^p$. 
This would imply that $F^p$, 
which maps each component of $\DD$ to itself,
maps $\gamma$ injectively to the unique path in $\BH$ that connects the endpoints of $\gamma$.
This would contradict Lemma \ref{no inv hor arcs}.
\end{proof}

We have the following two immediate corollaries:
\begin{cor}
   For any $\gamma\in \BH$, we have $\displaystyle{\sum_\de M^{pl}_{\gamma\de}\geq 2^l}$.

\end{cor}

\begin{cor} \label{cor:shortest-lp}
The shortest path in $\BH^{pl}$ connecting the endpoints of an arc in $\BH$ has length at least $2^l$. 
\end{cor}

\subsection{Lifting vertical arcs}

Let $\BH^\perp$ be the arc diagram consisting  of vertical  arcs in $\C\sm \DD$
that do not intersect the Hubbard tree $H$ (up to homotopy).
We observe that a horizontal arc on $\C \sm \DD$ lies in $\BH$ if and only if it does not intersect $\BH^\perp$. 

For any vertical arc $\beta$, 
we let $F^*\beta$ denote the set of lifts of $\beta$. 
Likewise for any nonempty set $\boldsymbol{\beta}$ of disjoint vertical arcs,
we let $F^*\boldsymbol{\beta}$ be the union of $F^*\beta$ over all $\beta \in \boldsymbol{\beta}$. 
We observe that $F^*\boldsymbol{\beta}$ will also be a nonempty set of disjoint vertical arcs. 
We observe that $F^*\BH^\perp = \BH^\perp$. 
We let $F^{l*}(\beta)$ be the result of $l$ iterations of $F^*$.

We will prove the following in this subsection. 
\begin{lemma} \label{lem:pull-to-get-Hperp}
For any vertical arc $\beta$, we have 
$$
\BH^\perp \subset \bigcup_{l \in \N} F^{*l}\beta. 
$$
\end{lemma}

Let us first define $|\beta|$ for any vertical arc $\beta$.
For any vertical path $b$,
we let $\underline b$ be the portion of $b$ from where it leaves from $\DD$ to the last place it intersects the actual Hubbard tree $H_F$. 
Then $\underline b$ is a bounded segment. 
We let $|\beta|$ be the infimum\footnote{This is finite by Lemma \ref{lem:un-wild}.}
of $|\underline b|$ over all representatives $b$ for $\beta$ (in $\C \sm \DD$).
We then observe the following.
\begin{enumerate}
\item
$|\beta|  = 0$ if and only if $\beta \in \BH^\perp$.
\item \label{it:contract}
If $\beta'$ is a lift of $\beta$, then $|\beta'| \le \lambda |\beta|$ for $\lambda < 1$ depending only on an upper bound for $|\beta|$.
\item \label{it:eta}
There is an $\eta > 0$ such that if $|\beta| < \eta$, then $\beta \in \BH^\perp$.
\end{enumerate}

We can now take our first step towards Lemma \ref{lem:pull-to-get-Hperp}. 
\begin{lem}\label{lem:eventually-H-perp}
For any vertical arc $\beta$ there exists an $n$ such that $F^{n*}\beta \subset \BH^\perp$. 
\end{lem}

\begin{proof}[Proof of Lemma \ref{lem:eventually-H-perp}]
Given $\beta$
take $\lambda = \lambda(|\beta|) < 1$ satisfying Property \ref{it:contract} above, and
let $n$ be the least integer greater than $\log_\lambda(\eta/|\beta|)$, where $\eta$ is given by Property \ref{it:eta}. 

Then if $\beta_0 = \beta, \beta_1, \ldots, \beta_n$ is a sequence of vertical arcs such that each $\beta_{i+1}$ is a lift of $\beta_i$,
we have $|\beta_n| \le \lambda^n |\beta| < \eta$, and hence $\beta_n \in \BH^\perp$. 
\end{proof}

We also observe
\begin{lemma} \label{lem:vert-to-vert}
If $\beta \in \BH^\perp$,
then $\bigcup_{k < 2p} F^{k*}\beta \supset \BH^\perp$.
\end{lemma}

\newcommand{\bhp}{\BH^\perp}
\begin{proof}
For each $k \in 0, \ldots, p-1, p$,
we let $\bhp_k$ be the vertical arcs in $\bhp$ that land on $D_k$ (with $D_p = D_0$).
We observe that $F^*\bhp_k = \bhp_{k-1}$ (for $k >0$),
that $|\bhp_k|$ is either 1 or 2 for each $k$, 
and that $|\bhp_k|$ is increasing in $k$ for $k \in [1, p]$,
and therefore $|\bhp_0| = 2$ and $|\bhp_1| = 1$.

Let $m < p$ be such that $\beta \in \bhp_m$.
We then have $F^{(m-1)*}\beta \subset \bhp_1$, 
and hence $F^{(m-1)*}\beta = \bhp_1$,
and therefore $F^{m*}\beta = \bhp_0= \bhp_p$, 
and $F^{m+p - k}\beta = \bhp_k$. 
\end{proof}

Putting together Lemmas \ref{lem:eventually-H-perp} and \ref{lem:vert-to-vert}, we immediately obtain Lemma \ref{lem:pull-to-get-Hperp}.

\section{Loss of horizontal weight} \label{sec:entropy} \label{sec:loss-hor}
In this section we prove the following purely combinatorial theorem,
which will be used to prove Theorem \ref{loss of hor weight} and ultimately Theorem \ref{def vertical teaser}. 
\begin{theorem}\label{loss-F}
Suppose $F(z) = z^2 + c$ is critically periodic of period $p$.
Let $\DD$ be the closure of the union of immediate attracting basins, and assume that $\DD$ is a union of disjoint closed Jordan domains.
Suppose that $X^n$, for $n \in \N$, are horizontal weighted arc diagrams on $\C \sm \DD$ such that
\begin{enumerate}
\item
$X^{n+1} \le X^n$, and 
\item
$F^*X^n \multimap X^{n+1}$
\end{enumerate}
for each $n \in \N$. 
Then
$$
     \| X^{8p}\|_1\leq \frac{1}{2} \| X^0 \|_1.
$$
\end{theorem}

Throughout this section,
we assume that the $\hx n$ (for $n = 0, \ldots, 8p$) are  horizontal weighted arc diagrams on $\C \sm \DD$
that satisfy Hypotheses 1 and 2 in Theorem \ref{loss-F}.

\subsection{Topological arrow}
  Let us consider two Riemann surfaces, $U\subset V$,
and let $\bal$ and $\bbe$ be multiarcs on $U$ and $V$  respectively.
We say that $\bal$ {\it topologically arrows} $\bbe$, $\bal\ta\bbe$,
if for any arc $\beta\in \bbe$
there is a sequence $(\alpha_k)$ of arcs
 with $\alpha_k \in \bal$ for each $k$,
and $(\alpha_k) \arrow \beta$.

A basic example comes from two WAD's, $X\in W(U)$ and $Y\in W(V)$,
such that $X\dominate Y$. Then  $\supp X \ta \supp Y$,
as follows immediately from the definitions.
So in the setting of this section, we have $$F^*\supp \hx {n} = \supp F^*\hx {n} \ta \supp \hx {n+1}$$ for each $n$. 

\subsection{An invariant horizontal arc diagram}

Let us say that a horizontal arc diagram $\bal$ on  $\C\sm \DD$ is {\it invariant} if
$$
   F^*\bal \ta \bal $$

\begin{prop}\label{invariance}
  There exists $n\leq 3p$ such that the horizontal arc diagram $\bal^n = \supp  \hx n$ is invariant.
\end{prop}

\begin{proof}
  Since $ \hx {n+1} \leq  \hx n$, we have
$\bal^{n+1}\subset \bal^n$.
Since $|\bal^0|\leq 3p$, there exists an  $n\leq 3p$ such that $\bal^{n+1}=\bal^n$. We take this $n$ as our desired $n$. 
Then, as we have already observed, we have $F^* \bal^n \ta \bal^{n+1} = \bal^n$.
\end{proof}

\subsection{Alignment with the Hubbard tree}
\newcommand{\bbb}{\boldsymbol \beta}
We begin with a lemma relating topological arrow (of horizontal arc diagrams) and vertical arcs. 
We observe that we can think of any vertical arc $\gamma$ in $\C \sm F^{-1}\DD$ landing at a disk in $\DD$ as a (unique and well-defined) vertical arc in $\C \sm \DD$.
\begin{lem} \label{lem:ta-int}
Suppose that 
\begin{itemize}
\item
$\bal$ is a horizontal arc diagram on $\C \sm F^{-1}\DD$,
\item
$\bbb$ is a horizontal arc diagram on $\C \sm \DD$, and 
\item
$\bal \ta \bbb$.
\end{itemize}
Suppose that $\gamma$ is a vertical arc on $\C \sm F^{-1}\DD$, landing at a disk in $\DD$, such that $\bal$ does not intersect $\gamma$.
Then $\bbb$ does not intersect $\gamma$.
\end{lem}
\begin{proof}
For every arc $\beta \in \bbb$, 
we can find a representing path $b$ such that the restriction of $b$ to $\C \sm F\inv\DD$ is a union of disjoint paths representing arcs in $\bal$. 
We can then realize $\gamma$ as a path (in $\C \sm F\inv\DD$)
 that does not intersect any of these $\bal$-representing paths. 
Therefore $\beta$ does not intersect $\gamma$.

We conclude that $\bbb$ does not intersect $\gamma$.
\end{proof}

\begin{lem}\label{H-perp}
  Any invariant horizontal arc diagram  $\bal$
is aligned with the Hubbard tree.
\end{lem}

\begin{proof}
Let  $\bal$ be an invariant horizontal arc diagram for a superattracting polynomial $F$,
and think of it for the moment as a union of actual paths,
each of which extends continuously to a map of a closed interval into $\C \sm \Int \DD$. 
Let $\De$ be the unbounded component of $\C\sm (\bal \cup \DD$).
If $\di \De$ were a subset of $\bal$,
then $\bal$ would contain a closed loop, which is absurd.
Therefore $\DD \cap \di\De$ is non-empty,
and we can find a path $b$ in $\De$ connecting $\infty$ to $\DD$. 
It represents a vertical arc $\beta_0$ such that  $\langle\bal, \beta_0\rangle=0$.

Now suppose that $\beta$ is any vertical arc for which $\langle\bal, \beta\rangle=0$,
and $\beta'$ is any lift of $\beta$ (that lands at a disk in $\DD$).
Then $\langle F^*\bal, \beta'\rangle=0$,
and, 
by Lemma \ref{lem:ta-int} and our assumption that $F^*\bal \ta \bal$,
we have $\langle \bal, \beta'\rangle=0$.
We conclude that $\pair\bal{F^*\beta} = 0$.

By induction,
we then have 
$\pair\bal{F^{*l}\beta} = 0$ for all $l$,
and hence, 
by Lemma \ref{lem:pull-to-get-Hperp},
$\pair\bal{\BH^\perp} = 0$. 
We conclude that $\bal \subset \BH$. 
\end{proof}
Putting this together with Proposition \ref{invariance} and Hypothesis 1 for the $\hx n$,
we obtain:

\begin{cor}\label{alignment}
  For any  $n\geq 3p$, the horizontal arc diagram $\bal^n = \supp  \hx n$ is
aligned with the Hubbard tree $H$.
\end{cor}

\subsection{The entropy argument}
We continue with the proof of Theorem \ref{loss-F},
where now we know that $\hx n$ is aligned with the Hubbard tree for each $n \ge 3p$. 
We first observe
\begin{lemma} \label{lem:dom-iter}
Under our standing assumption that $F^* \hx n \doms \hx {n+1}$,  we in fact have 
$F^{l*} \hx n \doms \hx {n+l}$ for all $n \in \N$ and $l \in Z^+$.
\end{lemma}
\begin{proof}
We proceed by induction,
assuming that $F^{l*} \hx n \doms \hx {n+l}$ and proving $F^{(l+1)*} \hx n \doms \hx {n+l + 1}$. 
We first observe that the assumption implies $$F^{(l+1)*} \hx n \doms F^* \hx {n+l},$$ 
and of course $$F^* \hx {n+l} \doms \hx {n+l+1}$$ by our standing assumption.
So we are finished by Lemma \ref{lem:dom-trans}.
\end{proof}

We then have, using the definition of $X|_D$ given in the preamble to Lemma \ref{lem:dom-restr},
\begin{lemma} \label{lem:res-iter}
Under our standing assumption,
$$\hx {n+l}|D_i \le (\deg F^l|D_i)\, \hx n|D_{i+l \bmod p}.$$
\end{lemma}
\begin{proof}
We have $F^{l*}\hx n|D_i = (\deg F^l|D_i)\, \hx n|D_{i+l \bmod p}$.
Moreover $F^{l*}\hx n \doms \hx {n+l}$ by Lemma \ref{lem:dom-iter}.
The Lemma then follows, by Lemma \ref{lem:dom-restr}.
\end{proof}
\begin{cor} \label{cor:sup-inf}
We have 
$$
\sup_i \hx {3n}|_{D_i} \le 2 \inf_i \hx {2n}|_{D_i} \le 2 \inf_i \hx 0|_{D_i}.
$$
\end{cor}
\begin{rem}
Lemma \ref{lem:res-iter} and its corollary is closely related to Theorem 9.3 of \cite{McM},
which relates the lengths of the hyperbolic geodesics around the small Julia sets of a given period. 
\end{rem}
We now have the following lemma, which is central to our entropy argument:
\begin{lemma} \label{lem:entropy-center}
For every $l \in \Z^+$, and $\beta \in \supp \hx {(3+l) p}$,
we have
$$
\hx {(3+l) p}(\beta) \le 2^{-l} \sup_i \hx {3p}|_{D_i}.
$$
\end{lemma}
\begin{proof}
We have $\supp X^{(3+l)p} \subseteq \BH^{lp}$ because $\supp X^{3p} \subseteq \BH$.

Suppose that $\beta$ connects $D_-$ and $D_+$.
By Lemma \ref{lem:into-flags},
there is a unique vertex-injective path in $\BH^{lp}$ connecting $D_-$ to $D_+$. 
We let the vertices of this path be $D'_0 = D_-, D'_1, \ldots, D'_n = D_+$,
and let the edges be $\gamma_0, \ldots, \gamma_{n-1}$,
where of course $D'_k \in \DD^{lp}$ and $\gamma_k \in \BH^{lp}$ for each $k$.
Moreover, 
by the same Lemma,
whenever $(\alpha_i) \to \beta$,
we can find, for each $k$, an $i$ such that $\alpha_i \in \Flag(D'_k, \gamma_k; \BH^{lp})$.
We also observe that $n \ge 2^{lp}$ by Corollary \ref{cor:shortest-lp}.

Choosing $k \in 0, \ldots, n-1$,
letting $j \in 0, \ldots, p-1$ be such that $F^{lp}(D'_k) = D_j$, and
letting $\eta$ be $F^{lp}(\gamma_k)$,
we observe that 
$$F^{lp}\from \Flag(D'_k, \gamma_k; \BH^{lp}) \to \Flag(D_j, \eta; \BH)$$
is a bijection.
Hence, 
$$
\norm{F^{lp*}\hx {3p}|_{\Flag(D'_k, \gamma_k; \BH^{lp})}}_1 = \norm{\hx {3p}|_{\Flag(D_j, \eta; \BH)}}_1
\le {\hx {3p}|_{D_i}} \le \sup_i {\hx {3p}|_{D_i}}.
$$

We then apply Lemma \ref{lem:in-each},
note the lower bound on $n$, and observe that $\bigoplus_{i=0}^{n-1} w_i \le w/n$ when $w_i \le w$ for each $i$.
\end{proof}

We can now prove Theorem \ref{loss-F}.
\begin{proof}
By Lemma \ref{lem:entropy-center} and Corollary \ref{cor:sup-inf},
letting $l = 5$, 
we have 
$$
\hx {8p}(\beta) \le  2^{-5} \sup_i \hx {3p}|_{D_i} \le 2^{-4}\inf_i \hx {0}|_{D_i}.
$$
The arc diagram $\supp \hx {8p} \cup \BH^\perp$ has at most $3p$ arcs and $\BH^\perp$ has at least $p$ arcs,
so $\hx{8p}$ has at most $2p$ arcs. 
We can then conclude
$$
\norm{\hx {8p}}_1 \le 2p \cdot 2^{-4} \inf_i \hx{0}|_{D_i} = \frac p 8 \inf_i \hx{0}|_{D_i}.
$$
On the other hand, 
$$
2\norm{\hx 0}_1 = \sum_i \hx 0|_{D_i} \ge p \inf_i \hx 0|_{D_i}. \qedhere
$$
\end{proof}

\section{Horizontal and vertical weight} \label{sec:restriction} \label{sec:hor-vert}
We now return to the setting of a quadratic-like map $f\from U'\to U$.
where $f$ is primitively renormalizable with some period $p$.
We let $\KK$ be the cycle of the corresponding little (filled) Julia sets.
We let $U^n = f^{-n}(U)$ (so $U^0 = U$ and $U^1 = U'$),
and observe that $U^{n+1} \subset U^n$,
and $f\from U^{n+1} \to U^n$ is quadratic-like.
In what follows we will relate the canonical WAD's of these surfaces.
For any $k, n \in \N$,
we will think of an arc in $U^n \sm f^{-k}\KK$ as horizontal if it connects two components of $f^{-k}\KK$,
and vertical if it connects one such component with $\di U^n$.

In this section we will prove Theorem \ref{def vertical teaser}, which relates the horizonal and vertical parts of $\Wcan(U\sm\KK)$. 
In Section \ref{domination sec}, we will use Theorem \ref{loss-F} to obtain Theorem \ref{loss of hor weight}, 
a loss of horizontal weight as we pass from $U \sm \KK$ to $U^{8p} \sm \KK$. 
In Section \ref{subsec:vert-from-hor-loss}, we will use this loss of horizontal weight to obtain definitely proportional vertical weight for $U^{8p} \sm \KK$.
In Section \ref{covering-application}, using a version of the Covering Lemma that we will state and prove in Section \ref{subsec:covering-lemma},
we will push forward the vertical weight for $U^{8p} \sm \KK$ back to $U \sm \KK$, to prove Theorem \ref{def vertical teaser}.
\subsection{The loss of horizontal weight in $W_\can$}\label{domination sec}
In this subsection we will produce a system of horizontal arcs that satisfy the hypothesis of Theorem \ref{loss-F}, 
and use that Theorem to prove Theorem \ref{loss of hor weight}. 

Let $\tl \KK=f^{-1}(\KK)$.
Then we have an embedding\footnote{The embedding here is of course the inclusion. In Sections \ref{pseudo-puzzle} and \ref{sec:conclusions} it will be replaced by an immersion, so it is good to have notation for it.}
 $i\from U^{n+1}\sm \KK \to U^n\sm \KK$,
a covering $f\from U^{n+1}\sm \tl \KK \to U^n\sm \KK$,
and an inclusion $U^{n+1}\sm \tl \KK \subseteq U^{n+1}\sm \KK$,
that together form a (non-commutative) triangle diagram.

\begin{center}
\begin{tikzcd}
 & U^n\sm \KK \\
U^{n+1}\sm \tl \KK  \arrow[r, "\subseteq"] \arrow[ur, "f"] & U^{n+1}\sm \KK \arrow[u, hook, "i"] \\
\end{tikzcd}
\end{center}

For $n \in \N$, 
we define $q_n\in \R_+$ letting $q_0 = 0$ and 
\begin{equation}\label{q}
                 q_{n+1}: =  3p(q_n+2).
\end{equation}
We observe that the $q_n$ are non-decreasing in $n$. 
Let $X^n = W_\can^\hor (U^n\sm \KK)$ and $\hat X^n= X^n- q_n$.

\begin{prop}\label{restrictions}
  We have:
\begin{itemize}

\item [(i)] $\hat X^{n+1} \leq i^* \hat X^n;$

\item [(ii)] $f^* \hat X^n\multimap \hat X^{n+1}.$

\end{itemize}
\end{prop}

\begin{proof} (i)
Since the embedding $i$ is proper on the ends corresponding to the little Julia sets $K_p(j)$,
we have  $X^{n+1} \leq i^* X^n$ by Property C.
This yields the desired inequality
since the sequence $(q_n)$ is non-decreasing.

\msk (ii)
Since the covering $f$ maps $U^{n+1}$ to $U^n$ and maps $\tl\KK$ to $\KK$,
horizontal and vertical  arcs of $U^n\sm \KK$ lift respectively to horizontal and vertical arcs
of $U^{n+1}\sm \tl \KK$. Hence  $f^* X^n = W_\can^\hor(U^{n+1}\sm \tl\KK)$.

By Property D,
$$
   W_\can(U^{n+1}\sm \tl\KK) \multimap W_\can^\hor (U^{n+1}\sm \KK)- 6p.
$$
But since the embedding  $U^{n+1}\sm \tl \KK \subset U^{n+1}\sm \KK$
is proper on $\di U^{n+1}$, the itinerary of any horizontal path $\gamma$ in $U^{n+1}\sm \KK$
consists only of horizontal arcs of $U^{n+1}\sm\tl\KK$.
It follows that
$$
  W_\can^\hor(U^{n+1} \sm \tl\KK) \multimap W_\can^\hor (U^{n+1}\sm \KK)-6p.
$$
Thus,
$
     f^* X^n \multimap X^{n+1} - 6p.
$
By Corollary \ref{cor:sub-sub} (and \eqref{q}), we then have
$$
  f^* \hat X^n = f^* X^n - q_n \multimap X^{n+1} -6p - 3p q_n \geq \hat X^{n+1}. \qedhere
$$
\end{proof}

We can now prove our fundamental lemma on the loss of horizontal weight under restriction:
\begin{theorem}\label{loss of hor weight}
  Let $f$ be a renormalizable quadratic-like map with period $p$.
Then there exists $M=M(p)$ such that
$$
     \| W^h_\can(U^{8p}\sm \KK) \|_1\leq \frac{1}{3} \| W^h_\can(U\sm \KK)\|_1,
$$
provided $ \| W^h_\can(U\sm \KK)\|_1>M(p)$.
\end{theorem}

\begin{proof}
As described in Section \ref{subsec:model},
we can transfer the $\hat X^n$ to horizontal weighted arc diagrams on $\C \sm \DD$ for the critically periodic model $F$ for 
$f$.
By Propositon \ref{restrictions},
these transferred $\hat X^n$ satisfy the hypotheses of Theorem \ref{loss-F},
so we can apply that theorem to obtain that $\norm{\hat X^{8p}}_1 \le \frac14\norm{\hat X^0}_1$, 
and therefore
\begin{align*}
   \| W^h_\can(U^{8p}\sm \KK) \|_1&\le \norm{\hat X^{8p}}_1 + 3pq_{8p} \\
&\le \frac14\norm{\hat X^0}_1 + 3pq_{8p} \\
&=    \frac{1}{4} \| W^h_\can(U\sm \KK)\|_1 + 3pq_{8p} \\
& \le \frac{1}{3} \| W^h_\can(U\sm \KK)\|_1
\end{align*}
when $ \| W^h_\can(U\sm \KK)\|_1 \ge 36p q_{8p}$. 
\end{proof}

\subsection{Vertical weight from the loss of horizontal} \label{subsec:vert-from-hor-loss}
We can now use Theorem \ref{loss of hor weight} to show that the vertical weight at level $8p$ is a definite proportion of the total weight at level 0. 

We first observe
\begin{lem}\label{increase of the total}
 For any $n\in \N$,
we have:
$$
  \| W_\can^{\ver+\hor}(U \sm \KK)\|_1 \leq \| W_\can^{\ver+\hor}(U^n \sm \KK)\|_1+ C(p).
$$
\end{lem}
\begin{proof}
This follows immediately from Theorem \ref{thm:total-increase},
because $U^n \sm \KK \subset U \sm \KK$,
and the inclusion is homotopic to a homeomorphism. 
\end{proof}

Combining Theorem \ref{loss of hor weight} with Lemma \ref{increase of the total}, we immediately obtain:
\begin{lemma} \label{separation property}
 Let $f$ be a renormalizable quadratic-like map with period $p$.
Then there exists $M=M(p)$ such that
$$
     \| W^\ver_\can(U^{8p}\sm \KK)\|_1 \geq  \frac{1}{2}  \| W_\can^{\ver+\hor}(U \sm \KK) \|_1 ,
$$
provided $ \| W_\can^h (U\sm \KK)\|_1>M(p)$.
\end{lemma}

\subsection{A suitable version of the Covering Lemma} \label{subsec:covering-lemma}
We will now state (and prove) a version of the Covering Lemma\cite{covering lemma} in the language of the canonical weighted arc diagram. In the next subsection, we will use Theorem \ref{thm:covering} to ``push forward'' the $\wcan^\ver(U^{8p} \sm \KK)$ 
to get a lower bound on $\wcan^\ver(U \sm \KK)$ and thereby prove Theorem \ref{def vertical teaser}.

As Theorem \ref{thm:covering} is simply a restatement of the Quasi-Invariance Law of \cite{covering lemma}, the reader may wish to read the statement of this Theorem and then skip to the next section. 
\begin{theorem} \label{thm:covering}
Suppose that 
\begin{itemize}
\item
$U$ and $V$ are simply connected Riemann surfaces,
\item
$A = \bigcup_{i=1}^p A_i$ is a union of $p$ disjoint compact connected full subsets of $U$,
\item
$B = \bigcup_{i=1}^n B_i$ is the same on $V$,
and 
\item
$f\from U\to V$ is a holomorphic branched cover of degree $D$,
\end{itemize}
such that 
\begin{itemize}
\item
all the critical values of $f$ lie in $B$,
\item
for each $i$,
$A_i$ is a component of $f^{-1}B_i$,
and 
\item
the degree of $f|_{A_i}$ is $d$
(in the sense that $A_i$ and $B_i$ have open neighborhoods $A_i'$ and $B_i'$ such that $f\from A_i'\to B_i'$ is a degree $d$ branched cover).
\end{itemize}
Suppose further that 
\begin{equation} \label{eq:covering-compare}
\norm{\wcan^\ver(U \sm A)}_1 \ge \epsilon \norm{\wcan^{\ver + \hor}(V \sm B)}_1,
\end{equation}
and 
\begin{equation} \label{eq:covering-large}
\norm{\wcan^\ver(U \sm A)}_1 \ge M(\epsilon, D, p)
\end{equation}
Then 
\begin{equation} \label{eq:covering-conclude}
\norm{\wcan^\ver(U \sm A)}_1 \le \frac{12d^2}{\epsilon}\norm{\wcan^\ver{V \sm A}}_1.
\end{equation}
\end{theorem}

Before proving this theorem,
let us review the notation of \cite{covering lemma} and relate it to $\Wcan(U \sm A)$.
Accordingly, 
let $U$ and $A$ be as in the statement of Theorem \ref{thm:covering}.
We then define three numbers determined by $U$ and $A$:
$$
  X = \WW(U,\, \bigcup_{j=1}^N A_j);
$$

\begin{equation}\label{Y}
   Y = \sum_{j=1}^N \WW(U, A_j),
\end{equation}

$$
   Z = \sum_{j=1}^N  \WW ( U \sm \bigcup_{k \not= j } A_k ,\, A_j ).
$$

It is easy to see that $X\leq Y\leq Z$.
We further claim that
\begin{equation} \label{eq:X-bounds}
\wcan^\ver(U \sm A) \le X \le \wcan^\ver(U \sm A) + 4p
\end{equation}
and
\begin{equation} \label{eq:Z-bounds}
\wcan^{2\hor + \ver}(U \sm A) \le Z \le  \wcan^{2\hor + \ver}(U \sm A) + 4p^2
\end{equation}
To see the lower bound for $X$ in \eqref{eq:X-bounds},
we need only observe that every vertical leaf in the canonical lamination for $U \sm A$ contributes to the defining width for $X$,
so $X$ is at least the total width of this vertical lamination.
To see the upper bound,
we consider the extremal foliation $\FF_X$ for the defining path family for $X$.
Then the leaves of $\FF_X$ can lie in at most $2p$ homotopy classes $\alpha_i$,
and 
$$
W(\FF_X|_{\alpha_i}) \le \wcan(U \sm A,\alpha_i) + 2
$$
by maximality. 
The derivation of \eqref{eq:Z-bounds} is similar and left as an exercise for the reader. 

We can now reduce Theorem \ref{thm:covering} to the version of the Quasi-Invariance Law stated in Section 3.1.1 of \cite{covering lemma}.
\begin{proof}[Proof of Theorem \ref{thm:covering}]
We wish to think of the setting of our Theorem as an instance of the setting of the QI Law stated in Section 3 of \cite{covering lemma},
and specialized in Section 3.1.1. 
So we apply the latter,
letting $U$ and $V$ in \cite{covering lemma} be our $U$ and $V$,
letting $\Lambda_j' = \Lambda_j$ be our $A_j$,
and letting $B_j' = B_j$ be our $B_j$.
Then the $Y_U$ in Section 3.1.1 of \cite{covering lemma} is our $Y(U, A)$; we will use only that it is at least $X(U, A)$.
\emph{Because the critical values of $f$ lie in $B$},
the $Z_V = Z_V\{B_j, B_j', CV\}$ in Section 3.1.1 is simply our $Z(V, B)$,
and of course $X_V$ in \cite{covering lemma} is our $X(V, B)$. 

The context and direction of the two inequalities in Section 3.1.1 of \cite{covering lemma} are such that we can replace $Y_U$ with $X_U$ in them. 
We can then derive,
assuming $\norm{W^\ver(U\sm A)} \ge 4p^\epsilon$,
\begin{align*}
Z_V &\le \norm{\wcan^{2\hor+\ver}(V \sm B)}_1 + 4p^2 \quad&\text{by \eqref{eq:Z-bounds}} \\
&\le 2 \norm{\wcan^{\hor+\ver}(V \sm B)}_1 + 4p^2 \\
&\le \frac3\epsilon \norm{\wcan^\ver(U \sm A)}_1 &\text{by \eqref{eq:covering-compare}}\\
&\le \frac3\epsilon X_U &\text{by \eqref{eq:X-bounds}}. 
\end{align*}
This is the separation assumption (with $\xi :=  3/\epsilon$) in Section 3.1.1 of \cite{covering lemma};
we can then conclude,
assuming \eqref{eq:covering-large} (which implies $X_U \ge K(\xi, p, D)$),
$$
X_U \le \frac{6d^2}\epsilon X_V,
$$
which implies
$$
\norm{\wcan^\ver(U \sm A)}_1 \le \frac{6d^2}\epsilon(\norm{\wcan^\ver(V \sm B)}_1 + 4 p).
$$
Assuming then that $\wcan^\ver(U \sm A) \ge 48d^2p/\epsilon$,
we obtain \eqref{eq:covering-conclude}.
\end{proof}

\subsection{Vertical part has a definite weight} \label{covering-application}
We can now prove Theorem \ref{def vertical teaser}:
\begin{proof}[Proof of Theorem \ref{def vertical teaser}]
   If $\| W_\can^\ver (U\sm \KK)\|_1 \geq \| W_\can^\hor (U\sm \KK)\|_1 $,
there is nothing to prove. So, we assume that the opposite inequality holds,
and then
$$
   \| W_\can^\hor (U\sm \KK)\|_1 > \frac{1}{2} \| W_\can^{\ver+\hor} (U\sm \KK)\|_1.
$$
Then taking $M=M(p)$ as in Corollary \ref{separation property},
we conclude that
\begin{equation}\label{separation for f^n}
     \| W^\ver_\can(U^{8p}\sm \KK)\|_1 \geq  \frac{1}{2}  \| W_\can^{\ver+\hor}(U \sm \KK) \|_1 ,
\end{equation}
provided  $\| W_\can (U\sm \KK)\|_1>2M$.

We are now ready to apply Theorem \ref{thm:covering} 
to the map $f^{8p}\from U^{8p} \to U$,
where we let $A = B = \KK$.
We observe that all the topological hypotheses are satisfied, 
with $d = 2^{8}$ and $D = 2^{8p}$, 
and also \eqref{eq:covering-compare} holds with $\epsilon = 1/2$.
Moreover,
taking $\normone{\wcan^{\hor+\ver}(U \sm \KK)}$ sufficiently large and applying Lemma \ref{increase of the total},
we have 
$$
\normone{\wcan^{\hor+\ver}(U^{8p} \sm \KK)} \ge \normone{\wcan^{\hor+\ver}(U \sm \KK)} - C(p) \ge M(\epsilon, D, p),
$$
as required in \eqref{eq:covering-large}.
We then obtain,
as our form of \eqref{eq:covering-conclude},
\begin{align*}
\normone{\wcan^{\ver}(U \sm \KK)} &\ge (24 \cdot 2^{16})^{-1} \normone{\wcan^{\ver}(U^{8p} \sm \KK)}.
\end{align*}
The Lemma then follows from this and Lemma \ref{separation property}.
\end{proof}

\section{Pseudo-quadratic-like maps and canonical renormalization} \label{pseudo-puzzle}
In this section we introduce the notion of a pseudo-polynomial-like map. 
This generalization of the polynomial-like-map is the result of the \emph{canonical renormalization}, 
which has a ``crucial property'' described in the Introduction and stated as Theorem \ref{thm:iterated-lengths};
this, 
along with an analog of Theorem \ref{thm:length-comparison} for pseudo-quadratic-like maps,
will allow us to prove the Main Theorem.

In Section \ref{pseudo:def} we define this new class of maps. 
In Section \ref{psi-pl} we show that every pseudo-polynomial-like map has a polynomial-like restriction. 
In Section \ref{pseudo-iterate} we describe how to iterate a $\psi$-polynomial-like-map. 
In Section \ref{can renorm} we define the canonical renormalization and state and prove Theorem \ref{thm:iterated-lengths}. 
The reader who finds the discussions in Section \ref{pseudo-iterate} and \ref{can renorm} too abstract and functorial for their tastes may prefer to skip to the end of each of these subsections, where there are more geometric perspectives. 
\subsection{Definition} \label{pseudo:def}
\begin{defn}
A \emph{pseudo-polynomial-like map}\footnote{Strictly speaking,
 we are defining a ``$\psi$-polynomial-like map with connected Julia set''.}
(or $\psi$-polynomial-like map)
of degree $d$ is a pair of holomorphic maps $(i, f)\from \BU'\to \BU$ and an FJ-set $K \subset \BU$ such that
\begin{enumerate}
\item 
$\BU'$ and $\BU$ are disks,
\item
$i$ is an immersion,
\item 
$f$ is a degree $d$ branched cover,
\item
$i^{-1}(K) = f^{-1}(K)$, and, denoting this set by $K'$,
\item
$K'$ is connected.
\end{enumerate}
\end{defn}
We will call $K$ the filled-Julia set of $(i, f)$. 
Note that $K'$ is connected if and only if
 all the branch values of $f$ lie in $K$.
Also note that,
because $f$ is proper,
$K'$ is an \fjset.
We will see that $K \subset \BU$ is uniquely determined by $(i, f)\from \BU' \to \BU$.
We will also see that these conditions imply that $i\from K'\to K$ is a bijection.

\begin{rem} \label{psi-plane}
By the Uniformization Theorem,
$\BU$ is isomorphic to either the complex plane or the unit disk.
In the former case, 
$\BU' = f^{-1}(\BU)$ is also the complex plane,
and $i\from \BU'\sm K'\to \BU\sm K$ is a holomorphic map from a punctured disk to a punctured disk, 
which then extends to the filled-in puncture by the Riemann removable singularity theorem.
It then follows that $i$ is an isomorphism, and then $f\circ i^{-1}$ is conjugate to a polynomial. 
\end{rem}

\subsection{Finding a polynomial-like restriction} \label{psi-pl}
We can now prove that every $\psi$-polynomial-like map $F$ has a polynomial-like restriction, 
with modulus controlled by the modulus for $F$. 
In the rest of this subsection we will assume that $d$ is an integer greater than 1,
and all constants and functions that appear will implicitly depend on $d$ (as well as any specified parameters). 

\begin{theorem}\label{ql extension}
Let $F=(i,f)\from\BU'\to \BU$ be a $\psi$-polynomial-like map of degree $d$
 with filled Julia set $K$.
Then  we can find $U \subset \BU$ open, with $K \subset U$
such that,
letting $U' \equiv f^{-1} U \subset \BU'$,
\begin{itemize}
\item
$f \from U'\to U$ is a degree d holomorphic cover,
\item
$i|_{U'}$ is an embedding,
\item
$i(U') \subset U$, and
\item
$\mod(U\sm i(U'))\geq \mu > 0$,
with $\mu$ uniformly positive for bounded $d$ and $\mod(\BU\sm K))$ bounded below.
\end{itemize}
\end{theorem}

\renewcommand{\aa}[1]{\A(1, #1)}
\newcommand{\aaa}[1]{\A(1/#1, #1)}

Before proving Theorem \ref{ql extension},
let us first prove a lemma about immersions of annuli.

\begin{lem} \label{lem:ann-immerse}
Let $d >1$ be an integer, take $r > 1$ and let $r' = r^{1/d}$. 
Suppose that $I\from \aa{r'} \to \aa{r}$ is a holomorphic immersion that extends continuously to $I\from \T \to \T$ of degree 1.
Take $\rho \in (1, r)$ and let $\rho' = \rho^{1/d}$.
\begin{enumerate}
\item \label{ai:s}
$I(\aa{\rho'}) \subset \aa s$,
where $s \equiv s(\rho, r)$ and $\rho - s$ is uniformly positive on any compact subset of $\{ (r, \rho) \mid 1 < \rho < r \le \infty \}$.
\item \label{ai:emb}
$I|_{\aa {\rho'}}$ is an embedding for $\rho \equiv \rho(r)$ uniformly greater than 1 when $r$ is bounded away from 1.
\end{enumerate}
\end{lem}

\begin{proof}
For Statement \ref{ai:s},
we can double $\aa {r'}$ and $\aa r$ over $\T$
(to form $\aaa {r'}$ and $\aaa{r}$), and extend $I$ to $I\from \aaa{r'}\to \aaa r$ by Schwarz reflection. 
Then the \poincare\ distance between $\T$ and $\T_{\rho'}$ in $\aa{r'}$
 is the same as that from $\T$ to $\T_\rho$ in $\aa r$. 
On the other hand, 
for $(\rho, r)$ in a compact subset of $1 < \rho < r < \infty$,
we must have uniform contraction of the \poincare\ metric, 
which proves our desired estimate in this case.
If $r = \infty$, 
then $I$ is the identity, and the estimate follows.
Finally, 
for $(\rho, r)$ in a compact subset of $1 < \rho < r \le \infty$,
we can imagine a sequence of such pairs $(\rho, r)$.
If the $r$'s are bounded, we are done by our first observation,
and if $r \to \infty$, we are done by the second.

For Statement \ref{ai:emb},
we observe that $I(\T_{\sqrt{r'}})$ is bounded away from $\T$ for $r$ bounded away from 1.
So,
for such $r$,
we have $\rho$ uniformly greater than 1 such that $I(\T_{\sqrt{r'}})$ is disjoint from $\Lambda := \aa \rho$.
Let $\Lambda'$ be the component of $I^{-1}(\Lambda)$ that is adjacent to $\T$;
then $I\from \T \cup \Lambda' \to \T \cup \Lambda$ is a covering map;
it must have degree 1, because $I\from \T \to \T$ does.
Therefore $I|_{\Lambda'}$ is an embedding,
Moreover,
by \poincare\ contraction on $\aaa {r'}$,
the outer boundary of $\Lambda'$ lies outside of $\aa{\rho'}$,
and therefore $I|_{\aa{\rho'}}$ is an embedding.
\end{proof}

\begin{rem} \label{ai:both}
Putting \eqref{ai:s} and \eqref{ai:emb} together,
we obtain, for $r$ bounded away from 1,
uniform $\rho$ and $s$ with $1 < s < \rho$ such that $I|_{\aa {\rho'}}$ is an embedding, 
and $I(\aa {\rho'}) \subset \aa s$. 
\end{rem}

\begin{rem} \label{ai:extends}
If we're given $I\from \aa{r'}\to \aa{r}$ such that $\abs{I(z)} \to 1$ as $\abs z \to 1$,
then $I$ extends continuously to $I\from \T \to \T$ by the Reflection Principle for harmonic functions;
see \cite{Ahlfors:complex} or \cite{Gamelin}.

If moreover $I\from \aa{r'} \to \aa{r}$ is injective in a neighborhood of $\T$,
then the extension to $\T$ is also injective, and hence has degree 1.
For if not,
we could find $z_0, z_1 \in \T$ with $z_0 \neq z_1$ and $I(z_0) = I(z_1)$,
and then,
by the Open Mapping Theorem, 
applied to the double of $I$ over $\T$,
every point in a neighborhood of $I(z_0)$ would be the image of a point near $z_0$ and a point near $z_1$. 
\end{rem}

\begin{proof}[Proof of Theorem \ref{ql extension}]
We first claim that there is a neighborhood of $K''$ on which $i$ is injective.
If not,
then there would be sequences $z_k \to K'$ and $w_k \to K'$ such that $i(z_k) = i(w_k)$,
and,
passing to subsequences such that $z_k \to z$ and $w_k \to w$,
we would have $i(z) = i(w)$.
This would imply that $z = w$ and that $z$ is a critical point of $i$, a contradiction.

We consider the annuli $\BU\sm K$ and $\BU'\sm K'$,
and uniformize them by the round annuli:
$$
\phi\from \A(1,r) \to \BU \sm K, \quad \phi'\from \A(1,r') \to \BU' \sm K',\quad \text{where}\ r'=r^{1/d}.
$$
Let $I= \phi^{-1}\circ i\circ \phi'\from \A(1,r')\to \A(1,r)$.
It is an immersion,
and,
by Remark \ref{ai:extends} and the previous paragraph,
$I$ extends continuously to $I\from \T\to\T$,
and the latter map has degree 1.

We can then find $1 < s < \rho$ as in Remark \ref{ai:both},
and let $U' = K' \cup \phi'(\aa {\rho'})$ and $U = K\cup \phi(\aa \rho)$. 
\end{proof}

\subsection{Iterating a $\psi$-polynomial-like map} \label{pseudo-iterate}
It is not immediately obvious how to iterate a $\psi$-polynomial-like map. 
We will first discuss the pullback in the category of complex 1-manifolds, 
and then use the pullback to define the iteration. 

\subsubsection{The pullback in the category of complex 1-manifolds}
Given two morphisms $g\from G\to V$ and $h\from H\to V$ in a given category,
the \emph{pullback} of these two morphisms is an object $E$ along with morphisms $\pg\from E \to G$ and $\ph\from E \to H$ 
such that $g \circ \pg = h \circ \ph$ 
and $(E, \pg, \ph)$ is universal in the sense that if $F$ is an object, and  $f_G\from F\to G$ and $f_H\from F\to H$ are morphisms
and $g\circ f_G = h\circ f_H$,  then there exists a unique morphism $\phi\from F\to E$
such that $f_G = \pg \circ \phi$ and $f_H = \ph \circ \phi$.
\begin{center}
\begin{tikzcd}
F \arrow[dr, "\phi" description,  dashed]  \arrow[drr, "f_h"] \arrow[ddr, "f_g",swap] & &\\
& E \arrow[r, "\ph"] \arrow[d, "\pg"] & H \arrow[d, "h"] \\
&G  \arrow[r, "g"] & V \\
\end{tikzcd}
\end{center}

This is of course old news for those who are at home in the world of  category theory, and invariably annoying for those who are not. It will help to consider two examples. In the category of sets, the pullback is just the equalizer: we let 
$$
E = \{ (x, y) \in G \times H \mid g(x) = h(y) \}
$$
and $\pi_G((x, y)) = x$ and $\pi_H((x, y)) = y$. 
In the category of topological spaces, we topologize $E$ as a subspace of $G \times H$ with the product topology. So the underlying set of the topological pullback is just the set-theoretic pullback. 

In any category, 
if the pullback exists (even just for a given diagram),
it is unique, 
in the sense that if $E$ and $E'$ are both pullbacks, there is an isomorphism $\psi\from E \to E'$ that makes everything commute. 
This follows immediately from the universal property of the pullback.

Let us observe and record a few properties of the topological pullback. 
The general principle is that a property of $g$ is inherited by $\ph$ (and likewise, by symmetry, a property of $h$ is inherited by $\pg$).
\begin{lemma} \label{lem:XYZ}
We have the following.
\begin{enumerate}
\item \label{it:proper}
If $g$ is proper, then $\ph$ is proper.
\item \label{it:loc-homeo}
If $g$ is a local homeomorphism, 
then $\ph$ is a local homeomorphism. 
\item \label{it:homeo}
If $g$ is a homeomorphism,
then $\ph$ is a homeomorphism. 
\item \label{it:set}
If $g$ is injective (or surjective, or bijective),
then $\ph$ is injective (or surjective, or bijective).
\item \label{it:set2}
For any $y \in H$, there is a natural 1-to-1 correspondence between $\ph^{-1}(\{y\})$ and $g^{-1}(\{ h(y) \})$. In particular, the two sets have the same cardinality. 
\end{enumerate}
\end{lemma}

\begin{proof} 
For Statement \eqref{it:proper},
suppose $B \subset H$ is compact.
Then
\begin{align*}
   \ph^{-1} (B) &= \{(x,y)\in G \times B \mid g(x) = h(y) \} \\
 &= \{ (x,y)\in g^{-1}(h(B)) \times B \mid g(x) = h(y) \} 
\end{align*}
which is a closed subset of a compact set.

For Statement \eqref{it:loc-homeo},
take $(x, y) \in E$.
Then we can find an open set $U \subset G$ such that $x \in U$ and $g|_U\from U \to g(U)$ is a homeomorphism.
Let $\hat U = \pg^{-1} U$. 
Then $\ph|_{\hat U}$ is a homeomorphism,
because its inverse is $y \mapsto ((g|U)^{-1}(h(y)), y)$. 

Statement \eqref{it:homeo} is essentially tautological, and holds in any category (where one of the two maps is an isomorphism).
Statements \eqref{it:set} and\eqref{it:set2} are straightforward (and belong to the category of sets) and are left to the reader. 
\end{proof}

We can now state and prove a lemma on a special case of pullback in the category of complex 1-manifolds. 
\begin{lemma} \label{lem:pullback}
The pullback of the diagram $g\from G\to V, h\from H\to V$,
in the category of complex 1-manifolds, 
exists when $g$ (or symmetrically, $h$) is an immersion. 
In this case the underlying topological space of the (Riemann surface) pullback is the topological pullback.
\end{lemma}
\begin{proof}
From Statement \eqref{it:loc-homeo} in Lemma \ref{lem:XYZ},
we have that $\ph\from E\to H$ is a local homeomorphism. 
We can therefore pull the complex 1-manifold structure (the system of charts) on $H$ back to $E$ (by $\ph$), 
and then $\ph$ is holomorphic. 

For any $(x, y) \in E$, 
we can proceed exactly as in the proof of Statement \eqref{it:loc-homeo} in Lemma \ref{lem:XYZ} 
to obtain $U \subset G$ and $\hat U \subset E$ with the stated properties. 
Then $\pg = g^{-1} \circ h \circ \ph$ on $\hat U$, as is therefore holomorphic there.
We conclude that $\pg$ is holomorphic on $E$.

It remains only to show the universal property.
Given $F, f_G, f_H$ as before,
there is a unique continuous map $\phi\from F \to E$ that makes the diagram above commute.
We need only show that $\phi$ is holomorphic,
and this follows from $\phi = (\ph|_{\hat U})^{-1} \circ f_H$ on $\phi\inv {\hat U}$ (where $\hat U$ is defined as before).
\end{proof}
We have written ``complex 1-manifold'' instead on ``Riemann surface" because the pullback, even in this special case, is not necessarily connected. For example, we can let $G$ and $H$ be disks in $V := \C$ with disconnected intersection; the pullback is then canonically isomorphic to $G \cap H$. 

Having proven this case let us observe the following:
\begin{lemma}
Suppose that $g\from G \to V$ and $h\from H\to V$ are holomorphic maps of complex 1-manifolds, 
and $g$ is an immersion and $h$  is proper of degree $d$.
Then in the pullback (with the notation as before), $\ph$ is an immersion, and $\pg$ is proper of degree $d$. 
Moreover the local degree of $\pg$ at each point $z$ in the pullback is equal to the local degree of $h$ at $\ph(z)$. 
\end{lemma}
\begin{proof}
We have already observed that $\ph$ is a local homeomorphism, and hence an immersion.
Likewise, by Statement \eqref{it:proper} of Lemma \ref{lem:XYZ}, we see that $\pg$ is proper. 
Moreover, choosing any point $(x, y) \in E$ such that $h(y)$ is not a critical value of $h$, 
we see by Statement \eqref{it:set2} (of Lemma \ref{lem:XYZ}) that $|\pg\inv {\{x\}}| = |h\inv {\{ h(y) \}}| = d$.
Therefore the degree of $\pg$ is $d$. 

The last statement then follows by considering small neighborhoods around $z$, $\pg(z)$, $\ph(z)$, and $g(\pg(z)) = h(\ph(z)$. 
\end{proof}

This is all that we will need in what follows, and reader may now skip to Section \ref{psiql:iterate-actual}.

Let us now discuss the general case of pullbacks for complex 1-manifolds,
as it may be of some interest and application.
We first consider an illustrative special case,
where $G = H = V = \D$, 
and $g(z) = z^{ma}$ and $h(z) = z^{mb}$, with $a$ and $b$ relatively prime.
In the case where $m$ is 1,
we can identify the equalizer $E$ with $\D$ by sending $(x, y)$ (with $x^a = y^b$) to the unique $z \in \D$ for which $z^b = x$ and $z^a = y$.
This of course gives us a Riemann surface structure on $E$ where the maps $\pg$ and $\ph$ are then $z \mapsto z^b$ and $z \mapsto z^a$ in our system of coordinates. 
To see universality, we take $F, f_G, f_H$ as specified above, and obtain a continuous map $\phi\from F \to E \equiv \D$ which we need to show is homeomorphic. 
This follows from Lemma \ref{lem:pullback} on $\phi\inv {\D \sm \{ 0 \}}$, and then follows on all of $F$ by the Riemann removable singularity theorem.

In the case where $m > 1$,
the topological (and set-theoretic) equalizer will be $m$ copies of $\D$ identified at their origins. 
We then need to desingularize this equalizer by separating these disks, 
and let our desired pullback be the desingularized equalizer.
So in this case, 
even though the topological pullback is connected,
the pullback in the category of complex 1-manifolds has $m$ components. 

More explicity we have $(x^a)^m = (y^b)^m$ 
and therefore $x^a = \rho^k_m y^b = (\rho^k_{mb} y)^b$, where $\rho_q := e^{2 \pi i/q}$, for some $k \in 0, \ldots m-1$. 
For each choice of $k$ we can then solve 
$$z^b = x \quad \text{and} \quad z^a = \rho^k_{mb} y$$
for a unique $z \in \D$.
So our complex 1-manifold pullback $\hat E$ is a disjoint union of $m$ disks parametrized by $k \in 0, \ldots m-1$;
we let $Z$ be the set of origins of these disks.
Given $F, f_G, f_H$ as before, 
we obtain a continuous $\phi\from F\to E$ as before.
Thinking of $E$ as the quotient of $\hat E$ where $Z$ is collapsed to a point $P$,
we see that $\phi$ maps each component\footnote{There is one component of $\phi^{-1}(E \sm P)$ for each component of $F$}
of $\phi^{-1}(E \sm P)$ to a component of $\hat E$. 
The Riemann removable singularity theorem then implies that $\phi$ has a unique continuous lift to $\hat \phi\from F \to \hat E$, 
and that this lift is in fact holomorphic.

We can now briefly sketch the proof of the general theorem.
\begin{theorem}
Pullbacks exists is the category of complex 1-manifolds.
\end{theorem}
\begin{proof}
For each point $(x, y)$ in the topological pullback, we can find coordinates around $x$, $y$, and $g(x) = h(y)$ such that $g$ and $h$ become $z\mapsto z^{ma}$ and $z \mapsto z^{mb}$ as in the preliminary discussion. We can then proceed as before to make the disingularized complex 1-manifold $\hat E$, and verify the universal property as in the previous discussion.
\end{proof}

\subsubsection{Iteration} \label{psiql:iterate-actual}
We now return to the setting of a degree $d$ $\psi$-polynomial-like map $(i, f)\from (\BU', K') \to (\BU, K)$,
with an eye towards constructing its iterates.
By Lemma \ref{lem:pullback},
we can form the pullback of the diagram $i\from \BU' \to \BU$, $f\from \BU' \to \BU$,
and we then label it as follows:

\begin{center}
\begin{tikzcd}
\BU'' \arrow[r, "f_1"] \arrow[d, "i_1"] & \BU'  \arrow[d, "i"]\\
\BU' \arrow[r, "f"]& \BU \\
\end{tikzcd}
\end{center}

By Lemma \ref{lem:XYZ},
$i_1$ is an immersion and $f_1$ is proper of degree $d$.
We can think of $\BU''$ as the set of $(x, y) \in \BU' \times \BU'$ such that $i(x) = f(y)$,
and then $i_1(x, y) = y$ and $f_1(x, y) = x$. 
Since 
$$x \in K' \iff i(x) \in K \iff f(y) \in K \iff y \in K',$$
we see that $i_1^{-1}(K') = f_1^{-1}(K')$; we call this set $K''$. 

To verify that $(i_1, f_1)\from (\BU'', K'') \to (\BU', K')$ is a degree $d$ $\psi$-polynomial-like map,
we now need only check that $\BU''$ is a disk and $K''$ is connected. 
Since $i\from K' \to K$ is a homeomorphism,
$i_1\from K'' \to K'$ is as well, and therefore $K''$ is connected.
By Lemma \ref{lem:XYZ},
$f_1\from \BU''\to \BU'$ is proper (of degree $d$),
and hence every component of $\BU''$ intersects $f^{-1}(K') = K''$. 
Therefore $\BU''$ is connected.
Moreover,
by the same Lemma,
$z \in \BU''$ is a critical point of $f_1$ if and only if $i_1(z)$ is a critical point of $f$, 
and the two have the same multiplicity.
Therefore all critical points of $f_1$ lie in $K''$,
and $i_1$ acts as a 1-to-1 correspondence between the critical points of $f_1$ and those of $f$, preserving multiplicity.
Hence the total number of critical points
(counting multiplicity)
of $f_1$ is $d-1$,
so,
by Riemann-Hurwitz,
the Euler characteristic of $\BU''$ is 1,
and hence $\BU''$ is a disk. 

We have thus proven
\begin{lemma} \label{lem:pullback-ps-plm}
The pullback $(i_1, f_1)\from (\BU'', K'') \to (\BU', K')$ 
of the degree $d$ $\psi$-polynomial-like map $(i, f)\from (\BU', K') \to (\BU, K)$ 
is itself a degree $d$  $\psi$-polynomial-like map.
\end{lemma}

Suppose that $(i, f)\from \BU' \to \BU$ is a $\psi$-polynomial-like map,
and $U' \subset \BU'$ and $U \subset \BU$ are such that $(i, f)\from U' \to U$ is a polynomial-like map,
in the sense that $i(U') \Subset U$ and $i|_{U'}$ is an embedding.
Letting $U''$ be the pullback of $(i, f)\from U' \to U$, 
we can think of $U''$ as a subset of $\BU''$,
and then $(i_1, f_1)\from U'' \to U'$ is a polynomial-like map.
The result is that the polynomial-like germ of $(i_1, f_1)\from \BU''\to \BU'$ is the same as that of $(i, f)\from \BU' \to \BU$. 

Let us relabel $(i, f)$ as $(i_0, f_0)\from (\BU^1, K^1) \to (\BU^0, K^0)$.
We can then apply Lemma \ref{lem:pullback-ps-plm} to construct degree $d$ $\psi$-polynomial-like maps
$$(i_n, f_n)\from (\BU^{n+1}, K^{n+1}) \to (\BU^n, K^n)$$ such that
$i_n \circ f_{n+1} = f_n \circ i_{n+1}$ for all $n$. 
By a very slight abuse of notation\footnote{or clever definition of $i$, $f$, and $K$ on the disjoint union $\bigcup_{n\in \N}\BU^n$},
we will often just denote $i_n$, $f_n$, and $K_n$ as $i$, $f$, and $K$.
This then allows us to refer to the iterations $i^n\from \BU^{m+n} \to \BU^{m}$ and $f^n\from \BU^{m+n} \to \BU^{m}$ (where we usually take $m=0$).

\subsubsection{Some other perspectives on iteration}
The pair $(i, f)$ provides a holomorphic map of $\BU'$ into $\BU \times \BU$, and the image of this map determines $\ff$.
So we can think of $\ff$ as this relation, and iterate it as we iterate any relation: the $n\thh$ iterate of a relation $R \subset X \times X$ is all pairs $(x_0, x_n)$ such that there exists $(x_0, x_1, \ldots x_n)$ for which $(x_i, x_{i+1}) \in R$ for each $i$. 
This iteration coincides with the iteration that we have defined through pullback.

A related perspective is to think of $i \circ f^{-1}$ as a map from $\BU$ to $\Sym_d \BU$, the symmetric $d$-fold product of $\BU$ with itself (or unordered $d$-tuples with multiplicity).
This map has no critical points, but it has $d-1$ ``anti-critical points'', where there image in $\Sym_d \BU$ has a point with multiplicity greater than 1. We can compose such maps by applying one map to each of the image points of the other, and taking the union. In this way the $n\thh$ iterate is a map from $\BU$ to $\Sym_{d^n}\BU$.

\newcommand{\hl}[1]{\hat #1}
\subsection{Canonical Renormalization}\label{can renorm}
We now describe a form of renormalization in the category of $\psi$-quadratic-like maps
that has the crucial property of preserving certain moduli and hyperbolic lengths. 

\subsubsection{Localization}
Let us first consider the category of \emph{marked Riemann surfaces}, 
which are pairs $(S, A)$, where $S$ is a Riemann surface,
and $A \subset S$ is compact and connected, 
and full in the sense no component of $S \sm A$ has compact closure in $S$. 
The morphisms of this category are holomorphic maps $g\from (S, A) \to (T, B)$ such that $A=g^{-1}(B)$. 

We say a marked Riemann surface is \emph{local} if every component of $S \sm A$ is an annulus. 
In this case $S$ deformation retracts into any neighborhood of $A$.
We observe that for every marked Riemann surface $(S, A)$ and neighborhood $V$ of $A$ there is a subsurface $U \subset V$ with $A \subset U$ such that $(U, A)$ is local. 
We will call such a $U$ a local neighborhood. 

A local marked Riemann surface $(S, A)$ for which $S \sm A$ is connected
is exactly a pair where $S$ is a disk and $A$ is an FJ-set.
Hence when $(i, f)\from (\BU', K') \to (\BU, K)$ is a $\psi$-polynomial-like map, 
then  $(\BU', K')$ and  $(\BU, K)$ are local marked Riemann surfaces, and $i$ and $f$ are morphisms. 
We will now define a functor from marked Riemann surfaces to local marked Riemann surfaces, 
that we will use in the next subsubsection to form the canonical renormalization.

We say that a morphism $g\from (S, A) \to (T, B)$
is a \emph{local cover} if $g\from A \to B$ is a bijection, and $g\from S \sm A \to T \sm B$ is a covering. 

\begin{theorem}
For any marked Riemann surface $(S, A)$, 
there is a unique (up to isomorphism) pair of a local marked Riemann surface $(S_A, \hat A)$ and local cover $\pi_{S,A}\from (S_A, \hat A) \to (S, A)$.

For any morphism $g\from (S, A) \to (T, B)$,
there is a unique continuous map $\hat g\from S_A \to T_B$ such that $\pi_{T, B} \circ \hat g = g \circ \pi_{T, A}$,
and $\hat g\from (S_A, \hat A) \to (T_B, \hat B)$ is a morphism. 

These operations define a functor from marked Riemann surfaces (and their morphisms) to the induced subcategory of local marked Riemann surfaces. 
\end{theorem}
\begin{proof}
For simplicity,
let us first consider the case where $S \sm A$ is connected;
this is the only case that we will actually use to make the canonical renormalization. 
To construct $S_A$ (and $\pi_{S, A}$),
we take a local neighborhood $U$ of $A$,
and take the cover $\pi'_{S, A}\from S'_A \to S \sm A$ that
corresponds to the fundamental group of $U \sm A$. 
Then there is a natural lift of $U \sm A$ to $S_A'$;
we let $S_A$ be the quotient of the disjoint union of $U$ and $S_A'$ by identifying points in $U \sm A$ with their image by this lift;
we let $\hat A \subset S_A$ be the image of $A$ by this quotient (and likewise define $\hat U \subset S_A$). 
Then $(S_A, \hat A)$ is a local marked Riemann surface, 
because $S_A \sm \hat A$ is isomorphic to $S'_A$,
which has fundamental group $\Z$.
We let $\pi_{S, A}(x) = \pi_{S, A}'(x)$ for $x \in S_A'$, and let $\pi_{S, A}$ be the inverse of the quotient map on $\hat U$.
These definitions are consistent, and make $\pi_{S, A}$ into a local cover.  

For the general case, 
we again take a local neighborhood $U$ of $A$ and
let $S'_A$ be the disjoint union of the annular covers of each component $E$ of $S \sm A$ induced by $E \cap U$.

Given a morphism $g\from (S, A)\to (T, B)$,
we let $\hat g = \pi_{T, B}|{\hat B}^{-1} \circ g \circ \pi_{S, A}$ on $A$,
and we let $\hat g$ be the lift of $g \circ \pi_{S, A}$ (by $\pi_{T, B}$) on every component of $S_A \sm \hat A$.
This defines $\hat g$ with $\pi_{T, B} \circ \hat g = g \circ \pi_{T, A}$,
and we can work backwards to see that $\hat g$ is the unique continuous map with this property.
One easily checks that $\hat g$ is a morphism.

The uniqueness of the commuting map $\hat g$ implies that lift of a composition of morphisms is the composition of the lifts; this shows functoriality. Finally, to see the uniqueness of the pair $(S_A, \hat A), \pi_{S, A}$, we suppose we are given two such pairs, and lift the identity map from $(S, A)$ to itself.
\end{proof}
We call  $(S_A, \hat A)$ the localization of $(S, A)$, 
and $\hat g$ the localization of $g$. 
We observe that the maps $\pi_{S, A}$ provide a natural transformation from localization 
(now thought of a functor from the category of marked Riemann surfaces to itself) to the identity functor.

Suppose $g\from (S, A) \to (T, B)$ is a morphism of marked Riemann surfaces,
and let $\hat g\from (S_A, \hat A) \to (T_B, \hat B)$ be its localization.
We observe that $\hat g$ is an immersion if $g$ is an immersion,
and that $\hat g$ is a branched cover (and $g$ is a branched cover) if $g\from S\sm A\to T \sm B$ is a covering map.

\newcommand{\vv}{\BU \sm B}
\newcommand{\vvp}{\BU' \sm f^{-1}(B')}
\subsubsection{Dynamical Localization}
Now let us suppose that 
\begin{itemize}
\item
$(i, f)\from (\BU', K') \to (\BU, K)$ is a $\psi$-polynomial-like map of degree $d$,
\item
$A \subset K$ is an \fjset\  (and let $A' = i^{-1}(A)$),
\item
$B \subset K$ is compact and $\BU \sm B$ is connected
\item
$A$ is a component of $B$
\item
$A'$ is a component of $f^{-1}(A)$
\item 
$i^{-1}(B') \subset f^{-1}(B)$, and
\item
all the branch values of $f$ lie in $B$.
\end{itemize}
Let $T = (\BU \sm B) \cup A$,
let $T' = (\BU' \sm f^{-1} B) \cup A'$.
Then $(T', A')$ and $(T, A)$ are marked Riemann surfaces,
and the restrictions of $i$ and $f$ to $T'$ are morphisms. 
So we can apply localization to obtain $(\hat i, \hat f)\from (T'_{A'}, \widehat{A'}) \to (T_A, \hat A)$,
and along with the local cover $\pi_{T, A}\from (T_A, \hat A)\to (T, A)$ (and the same for $(T', A')$). 
We let $\BV' = T'_{A'}$ and $\BV = T_A$. 
Since $A'$ and $A$ are FJ-sets and $(\BV', \widehat{A'})$ and $(\BV, \hat A)$ are local marked Riemann surfaces,
$\BV'$ and $\BV$ are disks.
Moreover $\hl i$ is an immersion
and $\hl f$ is a branched cover,
with degree equal to the local degree $d_{A'}$ of $f$ at $A'$.
So we have shown that $(\hl i, \hl f)\from (\BV', A') \to (\BV, A)$ is a $\psi$-polynomial-like map,
of degree $d_{A'}$. 
We call it the \emph{dynamical localization} of $f$ \emph{at} $A$ \emph{with the removal of} $B$. 
Moreover, 
since $\pi_{T,A} \from \BV \sm \hl A\to \BU \sm B$ is a covering,
the \poincare\ length of the core curve in $\BV \sm \hl A$
is equal to the \poincare\ length of the closed geodesic in $\BU \sm B$ that goes around $A$.

\subsubsection{Canonical Renormalization}
Now suppose that $(i, f)\from (\BU^1, K') \to (\BU^0, K)$ is a degree 2 $\psi$-polynomial-like map (which we will call a $\psi$-quadratic-like map).
Suppose that $f$ has a quadratic-like restriction that is primitively renormalizable of period $p$.
Then as before we let $K_p \subset K$ be the central small Julia set of the renormalization,
and $\KK_p$ be the union of small Julia sets.
We let $K_p' = i\inv{K_p}$. 
We observe that $K_p, \KK_p$ satisfy all the conditions for $A, B$ in the previous section,
when we replace $(i, f)$ with its iteration $(i^p, f^p)\from \BU^p \to \BU^0$. 
So we then let $(\hl i, \hl f)\from (\BV', \hl K_p') \to (\BV, \hl K_p)$ be the dynamical localization of $f^p$ at $K_p$ with the removal of $\KK_p$. 
We call this the \emph{canonical renormalization} of $f$ (or really $(i, f)$) of period $p$. 

As we did with polynomial-like maps,
if $(i, f)\from (\BU', K') \to (\BU, K)$ is a $\psi$-quadratic-like map,
then we let $\gamma(f)$ be the core curve of $\BU \sm K$.
Moreover, 
if $f$ is $p$-renormalizable,
we let $\gamma_p$ be the closed geodesic around $K_p$ in $\BU \sm \KK_p$. 
We can now prove the general version of the two ``crucial properties" that we stated in Section \ref{sec:intro}.
\begin{theorem} \label{thm:iterated-lengths}
Suppose that $(i, f)\from \BU' \to \BU$ is a $\psi$-quadratic-like map that is primitively renormalizable with periods $p$ and also $pq$ (with $p, q > 1$).
Let $(j, g)\from \BV' \to \BV$ be the canonical renormalization of $f$ of period $p$.
Then 
\begin{enumerate}
\item
$|\gamma(g)| = |\gamma_p(f)|$, and 
\item
$|\gamma_q(g)| \ge |\gamma_{pq}(f)|$.
\end{enumerate}
\end{theorem}
\begin{proof}
Statement (1) has already been observed in the definition of dynamical localization.

To prove (2),
we consider\footnote{We will write $K(g)$, $K_q(g)$, etc., for the objects associated with $g$, and let $K$, $\KK_p$, etc., be the objects for $f$.}
 the local cover 
$$\pi\from (\BV, K(g)) \to ((\BU \sm \KK_p) \cup K_p, K_p)$$
constructed for localization. 
We claim that $\pi(\BV\sm \KK_q(g)) \subseteq \BU \sm \KK_{pq}$;
this follows because $\pi(\BV \sm K(g)) \subseteq \BU \sm \KK_p \subseteq \BU \sm \KK_{pq}$
and $\pi(K(g) \sm \KK_q(g)) = K_p \sm \KK_{pq}$. 
Statement (2) then follows from the Schwarz Lemma.
\end{proof}

\subsubsection{Some other natural renormalizations} \label{satellite}

In the case where $f$ is twice renormalizable,
there is another natural renormalizations of $(i, f)\from \BU^1 \to \BU^0$ 
that we can form as dynamical localizations of $(i^p, f^p)\from \BU^p \to \BU^0$.
As before we let $A$ be $K_p \subset K$,
and let $B := A \cup \KK_{pq}$, where $f$ is also renormalizable of period $pq > p$.
We can do this even when $f$ is satellite renormalizable of period $p$;
this may be of value in the study of that case. 

There is a natural renormalization for which the $pq$-renormalization is equal to the $q$-renormalization of the $p$-renormalization (for a map that is $p$- and $pq$-renormalizable). Given $(i, f)\from \BU'\to \BU$, we let $P_f \subset \BU$ be the postcritical set. Then $\gamma_p$ (defined as before) separates $P_f \cap K_p$ from the rest of $P_f$. Choosing a base point lying inside of $\gamma_p$, we form the cover of $\BU \sm P_f$ associated to the subgroup of $\pi(\BU \sm P_f)$ represented by loops that never cross $\gamma_p$. We can then put $P_f \cap K_p$ ``back into'' this cover to form what we will call $\BV$; we can also think of $\BV$ as the disjoint union of this cover and the disk $D_p$ bounded by $\gamma_p$, quotiented by the natural embedding of $D_p \sm P_f$ into this cover. We observe that $K_p$ lifts to $\BV$; we will also call the lift $K_p$.

We can do the same with $f^{-p}P_f$ to form $\BV'$, and then $f^p$ lifts to a degree 2 proper map from $\BV' \to \BV$. Moreover the immersion of $\BU' \sm f^{-p}P_f$ into $\BU \sm P_f$ lifts to an immersion of $\BV' \sm (P_f \cap K_p)$ into $\BV \sm (P_f \cap K_p)$, and this induces an immersion of $\BV'$ into $\BV$. In this way we obtain a $\psi$-quadratic-like $p$-renormalization; one then verifies that the $pq$-renormalization is the $q$-renormalization of the $p$-renormalization.

\subsubsection{Paths and the canonical renormalization}
Suppose that $(i, f)\from \BU' \to \BU$ is primitively $p$-renormalizable and we want to produce the canonical renormalization $(j, g)\from \BV' \to \BV$.
For simplicity, let us assume that $K_p$ is locally connected. 
We can then let $\BV$ be the space of classes of paths that start on $K_p$ and end either on $K_p$ or in $\BV \sm \KK_p$, and in the former case can be deformed into $K_p$ in $V \sm (\KK_p \sm K_p)$. 
Two paths are equivalent if they are homotopic through paths with a common terminal point, and where the initial points stay in $K_p$. 
In $\BV$ there is a copy of $K_p$ (which we then think of as $K(g)$) given by the constant paths (which must map to a point in $K$).

We form $\BV'$ in the same way, with $(\BU^0 \equiv \BU, K_p, \KK_p)$ replaced by $(\BU^p, K_p, f^{-p}(\KK_p))$. We then observe that $j:=i^p$ applied to any path for $\BV'$ is a path for $\BV$, and equivalent paths are mapped to equivalent paths. Likewise for $g := f^p$. It is not hard to see that $(j, g)\from \BV' \to \BV$ is a $\psi$-quadratic-like map, and that this definition coincides with our definition in terms of localization.

\renewcommand{\BU}{\mathbf U}
\section{\emph{A priori} and beau bounds for bounded-primitive renormalization} \label{sec:conclusions}
We will begin this section by updating Theorems \ref{def vertical teaser} and \ref{thm:length-comparison} to the setting of $\psi$-quadratic-like maps.
We will then use the updated results to prove the \emph{a priori} and ``beau'' bounds for infinitely renormalizable maps of bounded-primitive type. 
This will be first expressed in the language of hyperbolic lengths of the geodesics around the small Julia sets,
and then in the form of the official modulus of quadratic-like maps that are formed as the renormalizations. 
\subsection{Estimates of extremal width and hyperbolic length for renormalizable $\psi$-quadratic-like maps}
Let us first prove the generalization of Theorem \ref{def vertical teaser} to $\psi$-quadratic-like maps:
\begin{theorem}\label{thm:pseudo-hor-ver}
Let $(i, f)\from \BU'\to \BU$ be a primitively renormalizable $\psi$-quadratic-like-map with period $p$
(and  $\KK$ is its cycle of little Julia sets).
If $\| W_\can^{\ver+\hor}(\BU\sm \KK)\|_1 > M(p)$ then
$$
   \|W^\ver_\can(\BU\sm\KK)\|_1 \geq C^{-1} \| W^\hor_\can(\BU\sm \KK)\|_1,
$$
where $M(p)$ depends only on $p$, and $C>0$ is an absolute constant.
\end{theorem}
\begin{proof}
We let $\BU^0 = \BU$, $\BU^1 = \BU'$, $f_0 = f$, and $i_0 = i$. 
As described in Section \ref{pseudo-iterate},
we can form disks $\BU^n$ and $\psi$-quadratic-like maps $(i, f)\from \BU^{n+1} \to \BU^n$.
Moreover these all have the same quadratic-like germ,
and are hence all $p$-renormalizable.
For each $n \in \N$,
we let $K^n \subset U^n$ be the central filled-Julia set of this map,
and let $\KK^n \equiv \KK_p^n \subset K^n$ be the union of small Julia sets (of the $p$-renormalization).
by a slight abuse of notation, 
we will just denote these by $K$ and $\KK \equiv \KK_p$.
We then have 
(with the slight abuse of notation described at the end of Section \ref{psi-pl}),
\begin{enumerate}
\item
a covering $f\from \BU^{n+1} \sm f^{-1} \KK \to \BU^n \sm \KK$,
\item
an inclusion $\BU^{n+1} \sm f^{-1} \KK \hookrightarrow \BU^{n+1} \sm \KK$, and
\item
an immersion $i\from\BU^{n+1} \sm \KK \to \BU^n \sm \KK$.
\end{enumerate}
The one significant change from the quadratic-like to the $\psi$-quadratic-like setting is the third map, 
which was previously an embedding, and is now an immersion. 

Now we need only go through Section \ref{sec:hor-vert} and see what needs to be changed.
In Section \ref{domination sec} we need only update the proof of Proposition \ref{restrictions} for the $\psi$-quadratic-like setting;
the proof of the updated Theorem \ref{loss of hor weight} is then identical. 
So we need only verify Properties (i) and (ii) of Proposition \ref{restrictions}.
Property (ii) follows just as before, 
because we use only the covering $f\from \BU^{n+1} \sm f^{-1} \KK \to \BU^n \sm \KK$ and the inclusion $\BU^{n+1} \sm f^{-1} \KK \hookrightarrow \BU^{n+1} \sm \KK$.
For Property (i) we need only observe that the immersion $i\from\BU^{n+1} \sm \KK \to \BU^n \sm \KK$ is still proper on the ends corresponding to $\KK$,
and Lemma \ref{e} applies to immersions as well as embeddings.

For Lemma \ref{increase of the total},
we can replace the inclusion $U^n \sm \KK \subset U \sm \KK$ 
with the immersion $i_n\from \BU^n \sm \KK \to \BU \sm \KK$;
this is still homotopic to a homeomorphism, 
so we can apply Theorem \ref{thm:total-increase} as before.
In a similar vein,
examining the proof of Theorem \ref{def vertical teaser} in Section \ref{covering-application},
we find that we need only replace $f^n\from U^n\to U$ with $f^n\from \BU^n \to \BU$,
where now this $f^n$ is the iteration described at the end of Section \ref{psi-pl}.
\end{proof}
We can now state and prove the version of Theorem \ref{thm:length-comparison} for a primitively $p$-renormalizable $\psi$-quadratic-like map $\ff \equiv (i, f)\from \BU' \to \BU$,
with $K$, $\KK$ defined as before, and geodesics $\gamma$ around $K$ and $\gamma_p$ around the central small Julia set in $\BU \sm \KK$:
\begin{theorem} \label{thm:pseudo-lengths}
Suppose that $\ff$ is a $\psi$-quadratic-like map that is primitively $p$-renormalizable.
Then, with the notation as described just above,
\begin{equation} \label{eq:pseudo-length-comparison}
|\gamma_p|> M(p) \implies |\gamma| \ge \epsilon p\, |\gamma_p|,
\end{equation}
where $M(p)$ depends only on $p$, and $\epsilon > 0$ is an absolute constant. 
\end{theorem}

\begin{proof}
This follows from Theorem \ref{thm:pseudo-hor-ver} in exactly\footnote{The proof of \eqref{eq:compare} (as Theorem 9.3 of \cite{McM}) can be readily adapted to $\psi$-quadratic-like maps.}
 the same way that Theorem \ref{thm:length-comparison} follows from Theorem \ref{def vertical teaser}. 
\end{proof}

The following is one of the main results of this paper. 
It is cited as Theorem 9.1 (as it appeared in the first arXiv version)
in \cite{decorations} and \cite{molecules}.

\begin{theorem}\label{improving of length}
For any $C \in \R^+$, there exists $\underline p$ such that for any $\bar p\geq \underline p$,
there exists $L(\bar p)$ with the following property.
Let $(i, f)\from \BU'\to \BU$ be primitively renormalizable $\psi$-quadratic-like map
with period $p$ such that  $\underline p\leq p\leq \bar p$.
Let $(j, g)\from \BV' \to \BV$ be the canonical renormalization of $f$ of period $p$.
Then (with the same notation as Theorem \ref{thm:iterated-lengths})
$$
   |\gamma(g)| > L(\bar p)\imply |\gamma(f)| > C |\gamma(g)|.
$$
\end{theorem}

We can now use Theorem \ref{thm:pseudo-lengths} to prove Theorem \ref{improving of length}.
\begin{proof}[Proof of Theorem \ref{improving of length}]
Given $C$, 
we let $\underline p = \lceil C/\epsilon \rceil$, 
where $\epsilon > 0$ is given by Theorem \ref{thm:pseudo-lengths}.
Then,
given $\bar p$,
we let $L(\bar p)$ be the maximum of $M(p)$ in Theorem \ref{thm:pseudo-lengths},
taken over all $p$ with $\underline p \le p \le \bar p$. 
Then
given the primitively $p$-renormalizable map $\ff$,
and letting $\gamma$ and $\gamma_p$ be as in Theorem \ref{eq:pseudo-length-comparison},
we have
\begin{align*} \label{eq:new-comparision}
|\gamma_p|> L(\bar p)& \implies 
|\gamma_p|> M(p) \implies 
|\gamma| \ge \epsilon p\, |\gamma_p|\\ 
& \implies
|\gamma| \ge \epsilon \underline p\, |\gamma_p| \implies
|\gamma| \ge C |\gamma_p|,
\end{align*}
where we have used \eqref{eq:pseudo-length-comparison} and our definitions of $L(\bar p)$ and $\underline p$.
Theorem \ref{improving of length} then follows from this and the first part of Theorem \ref{thm:iterated-lengths}.
\end{proof}

\subsection{\emph{A priori} and beau bounds for bounded-primitive renormalization}
We can now prove our ``if it's bad now, it was worse earlier'' statement for geodesics around the central small Julia sets.

We will have the following ``standard assumptions'' for this subsection.
We let $\ff$ be a $\psi$-quadratic-like map that is primitively renormalizable of periods $p_1\divides p_2\divides \ldots$,
with $1 < p_{k+1} / p_k \le B$ for all $k$.
The sequence of periods may be infinite,
or it may be finite,
in which case the results apply only when they refer to the $k$ for which $p_k$ exists. 
As explained in Section \ref{intro:defs},
we let $\gamma_{(k)}$ be the geodesic around $K_{p_k}(0)$ in $\C \sm \KK_{p_k}$;
we also let $\gamma_{(0)} = \gamma(\ff)$.
(If $\ff$ is a quadratic polynomial, then $\gamma_0$ does not exist, but we can let $| \gamma_0 | = 0$ and all results still apply.)

We can now prove a precise and general version of \eqref{eq:bad-now-worse-earlier}.
\begin{theorem} \label{thm:bnwe}
There exists $l$ such that for all $B$ there exists $M$:
Suppose that $\ff$ satisfies our standard assumptions.
Then,
for all $n$,
\begin{equation} \label{eq:bnwe}
|\gamma_{(n+l)}| > M \implies |\gamma_{(n)}| > 2|\gamma_{(n+l)}|.
\end{equation}
\end{theorem}

\begin{proof}
Choose $l$ such that $2^l \ge \underline p$, 
where $\underline p$ is given by Theorem \ref{improving of length} when we take $C := 2$.
Let $\bar p = B^l$,
and let $M = L(\bar p)$ in Theorem \ref{improving of length}. 
We claim that \eqref{eq:bnwe} holds for these choices (and all $n$).
To see this,
suppose that $|\gamma_{(n+l)}| > M$,
and let $\bgg \equiv (i, g)\from \BU' \to \BU$ be the canonical renormalization of $\ff$ of period $p_n$.
Let $q = p_{n+l}/p_n$;
then $\bgg$ is $q$-renormalizable,
and by Theorem \ref{thm:iterated-lengths},
we have 
\begin{equation} \label{eq:from-it}
|\gamma_{(n)}| = |\gamma(\bgg)| \text{ and } |\gamma_q(\bgg)| \ge |\gamma_{(n+l)}| > M.
\end{equation}
Since $\bar p \ge q \ge 2^l \ge \underline p$,
we can then apply Theorem \ref{improving of length} to $\bgg$ to obtain
$$
|\gamma(\bgg)| \ge 2 |\gamma_q(\bgg)|, 
$$
which combined with \eqref{eq:from-it} yields \eqref{eq:bnwe}.
\end{proof}

As a complement to Theorem \ref{thm:bnwe},
we observe that we can control the result of bounded-period renormalization. 
\begin{lemma} \label{lem:local-bound}
For all $P$ and $m > 0$ there is an $M \equiv M(m, P)$ such that 
if $\ff$ is a primitively renormalizable $\psi$-quadratic-like map of period $p \le P$ with $\gamma(\ff) \le m$,
then $|\gamma_p(\ff)| \le M$.
\end{lemma}
\begin{proof}
By Theorem \ref{thm:pseudo-lengths},
we have either $|\gamma_p| \le M(P)$ or $| \gamma_p | \le Pm /\epsilon$ for a universal $\epsilon$.
\end{proof}

Dennis Sullivan \cite{S} invented the term ``beau'', meaning \emph{bounded and eventually universally}\footnote{``beau'' looks and sounds a lot better than ``baeu''}. 
This is an extension and improvement on the \emph{a priori} bounds,
as the reader shall see in a moment. 
Here is the beau bounds for bounded-primitive renormalization, 
expressed in terms of the geodesics $\gamma_k$ around the central small Julia sets.
\begin{theorem} \label{thm:psi-beau}
For all $B$ there exists $M_0$ such that for all $m$ there exists $n_1$ and $M_1$:
Suppose $\ff$ satisfies our standard assumptions, and $| \gamma(\ff) | \le m$.
Then,
\begin{enumerate}
\item
for all $n\in \N$, $|\gamma_{(n)}| \le M_1$, and
\item
for $n \ge n_1$, $|\gamma_{(n)}| \le M_0$. 
\end{enumerate}
\end{theorem}

\newcommand{\mlong}{M_{\text{long}}}
\newcommand{\mshort}{M_{\text{short}}}

The first statement says that $| \gamma_{(n)} |$ is bounded in terms of $| \gamma(\ff) |$ (and $B$),
while the second statement says that $| \gamma_{(n)} |$ is eventually universally bounded (depending only on $B$),
with the number of renormalizations needed to become universally bounded depending on $| \gamma(\ff) |$. 

\begin{proof}
We take $l$ as in Theorem \ref{thm:bnwe},
and given $B$,
we let $\mlong$ be the $M(B)$ in Theorem \ref{thm:bnwe} 
and let $\mshort$ be $M(m, B^l)$ in Lemma \ref{lem:local-bound}.
($\mlong$ controls the long-term behavior, while $\mshort$ controls the short term.)

Let us first now Statement (1).
We take $M_1$ to be the greater value of $\mlong$ and $\mshort$.
Now suppose there were an $n$ with $|\gamma_{(n)}| > M_1$;
we take the least such value of $n$.
If $n \le l$ then $p_n \le B^l$,
and we contradict Lemma \ref{lem:local-bound}.
If $n > l$, 
then $|\gamma_{{(n-l)}}| \ge 2|\gamma_{(n)}|$ by Theorem \ref{thm:bnwe},
so $n$ was not minimal. 
This concludes the proof of (1).

For Statement (2) we take $M_0$ to be $\mlong$, which depends only on $B$.
We observe that we can rewrite Theorem \ref{thm:bnwe} as
\begin{equation} \label{eq:better}
| \gamma_{(n + l)} | \le \max(M_0, \frac12 | \gamma_{(n)} |).
\end{equation}
We let $n_1 = \ceil{1 + \log_2(M_0/\mshort)}\cdot l$. 
Now given $n \ge n_1$,
we can write $n = a + b l$,
where $a \le l$ and $b \ge \log_2(M_0/\mshort)$.
By Lemma \ref{lem:local-bound},
we have $| \gamma_{(a)} | \le \mshort$;
a simple induction with \eqref{eq:better} then implies that $| \gamma_{(a + kl)} | \le \max(M_0, 2^{-k} \mshort)$.
Letting $k = b$, we obtain our desired bound.
\end{proof}

\subsection{Quadratic polynomials and quadratic-like renormalization}
In this final subsection we return to the setting of quadratic polynomials and quadratic-like maps,
in order to derive the statements that will be more familiar (and sometimes more uto the reader. 

It will be convenient to define a \emph{semi-quadratic-like map} (with connected Julia set) to be a $\psi$-quadratic like map where the immersion is an embedding. In this case we can write it as $f\from U \to V$,
where $f$ is proper of degree 2 and $U \subset V$,
and the filled-Julia set $K \subset U$ is an FJ-set with $f^{-1}(K) = K$. 
Both quadratic-like maps and quadratic polynomials are semi-quadratic-like,
and indeed the main application of this terminology is to include both of these cases. 

In the case where $f$ is a quadratic polynomial,
the following is proven as Theorem 10.3 in \cite{McM} in the case where $M$ is sufficiently small,
and as Theorem 4.10 in \cite{towers} in general (for $f$ a with a slightly different statement that readily implies ours).
The proof here is fundamentally the same, but written with the machinery of the canonical renormalization.
\begin{theorem} \label{thm:gamma-to-modulus}
Suppose $f\from V' \to V$ is a primitively $p$-renomalizable semi-quadratic-like map,
and $|\gamma_p| \le M$.
Then $f$ has a (Douady-Hubbard) renormalization of period $p$ and modulus at least $m > 0$, for $m \equiv m(M)$. 
\end{theorem}
\begin{proof}
Let $V_p \subset V$ be the (open) disk bounded by $\gamma_p$;
and let $V'_p$ be $V_p$ pulled back along $\KK_p$,
so $V'_p$ is the unique component of $f^{-p}V_p$ that contains $K_p$.
Then $f^p\from V'_p \to V_p$ is a degree 2 cover;
we claim that $\cl{V'_p} \subset V_p$,
and that $\mod(V_p, V'_p) > m(M)$;
this would then complete the proof of the Theorem.
We let $A_p = V_p \sm K_p$ and $A'_p = V'_p \sm K_p$;
we need only show the outer boundary $\dout A'_p$  of $A'_p$ lies in $A_p$,
with a definite modulus between $\dout A'_p$ and $\dout A_p$.  

Now let $(i, g)\from \BU_p' \to \BU_p$ be the canonical renormalization of $f$ of period $p$.
We will also denote by $A_p$ and $A'_p$ their lifts to $\BU_p$ and $\BU'_p$ respectively. 
Then $g\from A'_p \to A_p$ is a degree 2 covering,
and $i\from A'_p \to A_p$ is an immersion that is degree 1 on the inner ideal boundary of $\di A'_p$
(see Remark \ref{ai:extends} and the proof of Theorem \ref{ql extension}). 
The inclusion $i(A'_p) \subset A_p$ (and lower bound on separating modulus) then follows from \eqref{ai:s} in Lemma \ref{lem:ann-immerse}.
\end{proof}
 
We can now state and prove the beau bounds for (semi-)quadratic-like maps and quadratic-like renormalization.
\begin{theorem} \label{thm:a-priori-bounds}
For all $B$ there exists $M_0$ such that for all $m$ there exists $n_1$ and $M_1$:
Suppose that $f\from U' \to U$ is a semi-quadratic-like-map that satisfies our standard assumptions, and $\abs{\mod(V \sm K)} \ge m$.
Then,
\begin{enumerate}
\item
for all $n\in \N$, 
we can find a quadratic-like $p_n$-renormalization $f^{p_n}\from V' \to V$ with $\mod(V\sm V') \ge M_1$, and 
\item
for $n \ge n_1$, 
we can find a quadratic-like $p_n$-renormalization $f^{p_n}\from V' \to V$ with $\mod(V\sm V') \ge M_0$.
\end{enumerate}
\end{theorem}
\begin{proof}
This follows immediately from Theorems \ref{thm:psi-beau} and \ref{thm:gamma-to-modulus}.
\end{proof}
Applied to quadratic-like maps, this is the standard statement of the beau bounds.

We conclude with a complete and explicit statement of the \emph{a priori} bounds,
which is an immediate corollary of Statement (1) of Theorem \ref{thm:a-priori-bounds}.
\begin{theorem}
For all $B$ there exists $m > 0$:
Let $f\from \C\to \C$ be an infinitely $B$-bounded primitively renormalizable quadratic polynomial,
with primitive periods $p_1 \divides p_2 \ldots$. 
Then for all $n$, there is a quadratic-like $p_n$-renormalization $f^{p_n}\from V' \to V$ with $\mod(V \sm V') \ge m$. 
\end{theorem}

\section{Appendix A: Extremal length and width}\label{apppendix: width}
[If you have come to this section from  \cite{decorations} or \cite{molecules},
looking for Theorem 9.1 of this paper,
it is now known as Theorem \ref{improving of length}.]

There is a wealth of sources containing background material on extremal length,
see, e.g., the book of Ahlfors \cite{A}.
We will briefly summarize the necessary minimum
(see also the Appendix of \cite{covering lemma}).

\subsection{Definitions} \label{A:def}
In this section we will think of a path as a smooth map $\gamma\from I \to U$,
where $I \subset R$ is an interval (that may be open or closed),
and $\gamma$ is not necessarily proper.

Let $\GG$ be a family of paths on a Riemann surface $U$.
Given a (measurable) conformal metric $\mu= \mu(z) |dz|$ on $U$,
let
$$   \mu (\GG) = \inf_{\gamma\in \GG}  \mu (\gamma).$$
 where $\mu(\gamma) $ stands for the $\mu$-length of $\gamma$.
The length of $\GG$  with respect to $\mu$ is defined as
$$
   \LL_\mu (\GG) = \frac {\mu(\GG)^2 } {\area_\mu (U)},
$$
where $\area_\mu$ is an area form of $\mu$.
Taking the supremum over all conformal metrics $\mu$, we obtain the {\it extremal length}
$\LL(\GG)$ of the family $\GG$.

  The {\it extremal width} is the inverse of the extremal length:
$$
 \WW(\GG)=\LL^{-1}(\GG).
$$
It can be also defined as follows. Consider all conformal  metrics $\mu$ such that
$\mu(\gamma)\geq 1$ for any $\gamma\in \GG$.  Then $\WW(\GG)$ is the infimum of
the areas $\area_\mu(U)$ of  all such metrics.

\subsection{Electric circuit laws}\label{electric}
We say that a family  $\GG$ of paths {\it overflows} a family $\GG'$ if
any path of $\GG$ contains some path of $\GG$.
Let us say that $\GG$ {\it disjointly overflows}  two families, $\GG_1$ and $\GG_2$,
if any path of $\gamma\in \GG$ contains the disjoint union $\gamma_1\sqcup \gamma_2$
of two paths $\gamma_i\in \HH_i$.

The following crucial properties of the extremal length and width show
that the former behaves like the resistance in electric circuits,
while the latter behaves like conductance.

\begin{serieslaw}
Let $\GG$ be a family of paths that disjointly overflows two other families,
                      $\GG_1$ and $\GG_2$.
Then
$$  \LL(\GG) \geq \LL(\GG_1) + \LL(\GG_2), $$
or equivalently,
$$  \WW(\GG)\leq \WW(\GG_1)\bigoplus \WW(\GG_2). $$
\end{serieslaw}

\begin{parallellaw}
  For any two families $\GG_1$ and $\GG_2$ of paths we have:
$$
    \WW(\GG_1\cup \GG_2)\leq \WW(\GG_1) + \WW(\GG_2).
$$
If $\GG_1$ and $\GG_2$ are contained in two disjoint open sets, then
$$
  \WW(\GG_1\cup \GG_2) = \WW(\GG_1) + \WW(\GG_2)
$$
\end{parallellaw}

\subsection{Transformation rules}\label{trans rules sec}
 Both extremal length and extremal width  are conformal invariants.
 More generally, we have:
\begin{lem}\label{increase}
   Let $f\from U\to V$ be a holomorphic map between two Riemann surfaces,
and let $\GG$ be a family of paths on $U$. Then
$$
      \LL(f(\GG))\geq \LL(\GG).
$$
\end{lem}
See Lemma 4.1 of \cite{covering lemma} for a proof.
\begin{rem}
We push forward \emph{parametrized} paths, so each $\gamma \in \GG$ is mapped to $f \circ \gamma \in f(\GG)$. 
For example, if $f\from U\to V$ is a degree 2 covering map between annuli (of finite modulus), and $\GG$ is the family of closed paths that go around $U$,
then $f(\GG)$ includes the closed paths that go twice around another closed path that goes once around $V$. 
While the extremal length of the closed paths that go once around $V$ is $\LL(\GG)/2$,
we actually have $\LL(f(\GG)) = 2 \LL(\GG)$ in this example.
\end{rem}

\begin{cor}\label{two families}
     Under the circumstances of the previous lemma,
let $\HH$ be a family of paths in $V$ satisfying the following lifting
property: any path $\gamma\in \HH$ contains an arc that lifts to some path
in $\GG$. Then $\LL(\HH)\geq \LL(\GG)$.
\end{cor}

\begin{proof}
   The lifting property means that the family $\HH$
overflows the family $f(\GG)$. Hence $\LL(\HH)\geq \LL(f(\GG))$,
and the conclusion follows.
\end{proof}

\begin{cor}\label{Q}
  Let $Q$ and $Q'$ be two quadrilaterals,
and let $e\from Q\to Q'$ be a holomorphic map that maps the horizontal sides of $Q$ to the horizontal sides of $Q'$.
Then $\WW(Q)\leq \WW(Q')$.
\end{cor}

\begin{proof}
  Let $\GG$ (resp. $\GG'$) be the family of horizontal paths in $Q$ (resp., in $Q'$).
 Since the horizontal sides of $Q$ are mapped to the horizontal sides of $Q'$,
these families satisfy the lifting property of the previous Corollary.
Hence $\LL(\GG)\leq \LL(\GG')$, and we are done.
\end{proof}

\begin{lem}\label{modulus transform-2}
 Let $f\from U\to V$ be a branched covering between two compact Riemann surfaces with boundary.
Let $A$ be an archipelago in $U$,  $B=f(A)$, and
assume that  $f\from  A\to B$ is a branched covering of degree $d$. Then
$$
      \mod(V, B) \geq d\, \mod (U,  A).
$$
\end{lem}

See Lemma 4.3 of \cite{covering lemma} for a proof.

Given two compact subsets $A$ and $B$ in a Riemann surface $S$,
let $\WW_S (A,B)$ stand for the {\it extremal width}
between them, i.e., the extremal width of the family of paths connecting $A$ to $B$.

\begin{lem}\label{modulus transform-1}
Let $S$ and $S'$ be two compact Riemann surfaces with boundary.
and let  $f\from S\to S'$ be a branched covering of degree $D$.
Let $S'=A'\sqcup B'$, where $A'$ and $B'$ are closed,
and let  $A=f^{-1}(A')$, $B=f^{-1}(B')$.
Then
$$
      \WW_S(A, B)= D \, \WW (A', B').
$$
\end{lem}

See \cite{A} for a proof.
It makes use of the fact that
the extremal width $\WW(A,B)$ is achieved on the {\it harmonic foliation} $\FF=\FF_S(A,B)$ connecting $A$ and $B$,
i.e.,  the gradient foliation of the harmonic function $\om$ vanishing on $0$ and equal to $1$ on $B$.
Hence $\WW_S(A,B)$ is equal to the $l_1$-norm of the associated  WAD $W_\FF$.

\subsection{Non-Intersection Principle}\label{so principle}

The following important principle says that
two wide quadrilaterals cannot go non-trivially one across the other:

\begin{lem}\label{no crossing}
Let us consider two quadrilaterals, $Q_1$ and $Q_2$,
endowed with the vertical foliations $\FF_1$ and $\FF_2$.
If $\WW(Q_i)\geq 1$ then \\
$\bullet$ either there exists a pair of disjoint leaves $\gamma_i$ of the foliations $\FF_i$;\\
$\bullet$ or $\WW(Q_1)=\WW(Q_2)=1$, and the rectangles perfectly match in the sense that
   the vertical sides of one of them coincide with the horizontal sides of the other.
\end{lem}

\begin{proof}
  Assume that the first option is violated,
so that, every leaf of $\FF_1$ crosses every leaf of $\FF_2$.
  Let us uniformize our rectangles by standard rectangles,
$\phi_i\from \BQ(a_i)\to Q_i$, where $a_i\geq 1$.
Let $\la_i$ be the Euclidean metrics on the $\BQ(a_i)$,
and let $\mu_i= (\phi_i)_*(\la_i)$.
Since every vertical leave $\gamma$ of $Q_2$ crosses every vertical leaf of $Q_1$,
$\mu_1(\gamma)\geq a_1$. Hence $a_2= \WW (Q_2)\leq 1/a_1$,
which implies that $a_1=a_2=1$. Thus, both $\BQ(a_i)\equiv \BQ$ are the squares.

Moreover, $\mu_1$ must be the extremal metric for $\FF_2$.
Since the extremal metric is unique (up to scaling), we conclude that
$\mu_1=\mu_2$. Hence $\phi_2^{-1}\circ \phi_1\from \BQ \to \BQ $
is the isometry of  the square, and the conclusion follows.
\end{proof}

\section{Appendix B: Six or seven lemmas about domination} \label{app:dom}
\subsection{Interchanging the order of arithmetic and harmonic sum}
We begin with an extremal problem that determines the harmonic sum:
\begin{lemma}
We have 
\begin{equation} \label{eq:energy-min}
\bigoplus_i v_i = \inf_{(w_i)} \frac{\sum_i w_i^2 v_i}{\left(\sum_i w_i\right)^2}
\end{equation}
where the infimum is taken over all sequences $(w_i)$ of non-negative numbers, not all zero. 
\end{lemma}
\begin{proof}
We first observe that we have equality when we let $w_i = 1/v_i$. Moreover, by the Cauchy-Schwarz inequality,
for any $(w_i)$, we have
$$
\left(\sum_i w_i\right)^2 \le \sum_i \frac 1{v_i} \sum_i w_i^2 v_i,
$$
and therefore
$$
\bigoplus_i v_i \le \frac{\sum_i w_i^2 v_i}{\left(\sum_i w_i\right)^2}. \qedhere
$$
\end{proof}
We can also write \eqref{eq:energy-min} as 
\begin{equation} \label{eq:energy-min-normalized}
\bigoplus_i v_i = \inf_{\sum_i w_i = 1} \sum_i w_i^2 v_i.
\end{equation}

We can now prove an inequality on interchanging the order of harmonic and arithmetic sum. 
\begin{lemma} \label{lem:sum-hsum}
We have
\begin{equation} \label{eq:sum-hsum}
\sum_i \bigoplus_j v_{ij} \le \bigoplus_j \sum_i v_{ij}.
\end{equation}
\end{lemma}
\begin{proof}
In light of \eqref{eq:energy-min-normalized},
we can reduce \eqref{eq:sum-hsum}  to 
\begin{equation}
\inf_{\forall i, \sum_j w_{ij} = 1} \sum_{ij} w_{ij}^2 v_{ij} \le \inf_{\sum_j x_j = 1} \sum_{ij} x_j^2 v_{ij},
\end{equation}
which is immediate, as given any $(x_j)$, we can let $w_{ij} = x_j$ for all $i$ and $j$. 
\end{proof}

\subsection{Two fundamental properties of domination}
If $X_1$ and $X_2$ are two weighted arc diagrams (on a surface $S$),
we let $\norm{X_1 - X_2}_\infty$ be the supremum of $|X_1(\alpha) - X_2(\alpha)|$ over all $\alpha$.
We can analogously define $\norm{X_1- X_2}_1$. 
We observe that $\norm{X_1 - X_2}_\infty \le \norm{X_1 - X_2}_1 \le 6|\chi(S)| \norm{X_1 - X_2}_\infty$,
so the topology on $\WW(S)$ defined by the two norms is the same.
We call it the strong topology on $\WW(S)$, 
where the weak topology is given by the product topology on $[0, \infty)^{\AA(S)} \supset \WW(S)$. 
We observe that $X_n \to X$ strongly implies $X_n \to X$ weakly.
We also observe that the converse holds if there is a finite set $E$ of arcs such that $\supp X_n \subset E$ for all $n$. 

Here now is our first fundamental lemma:
\begin{lemma} \label{lem:domination-closed}
Suppose $X_n \doms Y_n$ and $X_n\to X$ strongly and $Y_n \to Y$ weakly. Then $X \doms Y$. 
\end{lemma}

This follows almost immediately from the following more simply stated lemma.
\begin{lemma} \label{lem:y-minus-epsilon}
Suppose $X \doms Y-\epsilon$ for all $\epsilon > 0$. Then $X \doms Y$. 
\end{lemma}

\begin{proof}[Proof of Lemma \ref{lem:y-minus-epsilon}]
If $(\alpha_j)$ and $(\alpha'_k)$ are two sequences of arcs, 
we say that $(\alpha_j) \prec (\alpha'_k)$ if $\sum_j \alpha_j \le \sum_k \alpha'_k$ and the two sums are not equal.
We say that $(\alpha_j)_j$ \emph{efficiently arrows} $\beta$ if $(\alpha_j) \to \beta$ and there is no $(\alpha'_k)$ such that $(\alpha'_k) \prec (\alpha_j)$ and $(\alpha'_k) \to \beta$.
We can define the relation ``$X$ efficiently dominates $Y$'' by replacing $(\alpha_j) \to \beta$ with ``$(\alpha_j)$ efficiently arrows $\beta$'' in the definition of domination. 
We then observe that $X \doms Y$ actually implies that $X$ efficiently dominates $Y$,
because we can replace each summand in the domination relation with one that comes from the efficient domination. 

We now claim that for every $\beta \in \supp Y$, 
there are only finitely many sequences $(\alpha_j)$ (taken from the support of $X$) that efficiently arrow $\beta$.
To see this, let us consider the set of sums $\sum_j \alpha_j$ taken over all sequences $(\alpha_j)$ that arrow $\beta$.
Each such sum is an element of the free abelian monoid over the support of $X$, which is of course isomorphic to $\N^m$ for some $m$. 
We claim that \emph{any} subset of $\N^m$ has only finitely many minimal elements, which of course would imply the former claim.

To see this latter claim (which is a well-known property of $\N^m$) for a subset $T$ of $\N^m$,
we first observe that,
by induction on $m$,
there are,
for each $i \in 1,\ldots, m$ and $k \in \N$,
only finitely many minimal elements $(n_1, \ldots, n_m)$ of $T$ for which $n_i =k$.
Now take any $(n_1, \ldots, n_m)$ of $T$;
any minimal element $(n'_1, \ldots, n'_m)$ must have $n'_i = k$ for some $i \in 1, \ldots, m$ and $k \le n_i$. 

So now suppose that $X \doms Z$, and hence that $X$ \emph{efficiently dominates} $Z$.
Then, in light of the finiteness of efficiently arrowing sequences $(\alpha_j)$ (for each $\beta \in \supp Z$),
and applying Lemma \ref{lem:sum-hsum},
we can write $X \ge \sum X_{\ba_i, \beta_i}$ and $Y = \sum v_i \beta_i$,
where $X_{\ba_i, \beta_i} = \sum_j w_{ij} \ba_i(j)$ and $\ba_i = (\ba_i(j))_j$ efficiently arrows $\beta_i$ (and $\oplus_j w_{ij} \ge v_i$),
and  $(\ba_i)$ is a fixed finite sequence of sequences of arcs, only depending on $\supp Z$.

So now suppose that $X \doms Y -\epsilon$ for all $\epsilon$.
Then, for each $\epsilon$, we can find $X^\epsilon_{\ba_i, \beta_i}$ that organize the efficient domination of $X$ over $Y - \epsilon$.
We then pass to a subsequence of $\epsilon$ such that each of these terms converge. 
\end{proof}

Our next lemma states that domination is transitive.
\begin{lemma} \label{lem:dom-trans}
Suppose $X \doms Y$ and $Y \doms Z$. 
Then $X \doms Z$.
\end{lemma}
Before proving Lemma \ref{lem:dom-trans},
we will need a simple sublemma, which is left to the reader. 
\begin{lemma} \label{lem:dom-split}
Suppose that $X \doms Y$ and $Y = \sum Y_k$.
Then we can find $X_k$ such that $X \ge \sum X_k$ and $X_k \doms Y_k$ for each $k$.
\end{lemma}
\begin{proof}[Proof of Lemma \ref{lem:dom-trans}]
By Lemma \ref{lem:dom-split},
we can reduce to the case where 
\begin{enumerate}
\item
$Z = w \gamma$,
\item
$Y = \sum v_i \beta_i$,
\item
$X = \sum u_{ij} \alpha_{ij}$,
\item
$(\beta_i) \to \gamma$ and $\oplus_i v_i \ge w$, and
\item
for each $i$, $(\alpha_{ij})_j \to \beta_i$ and $\oplus_j  u_{ij} \ge v_i$. 
\end{enumerate}
We then have
$(\alpha_{ij})_{(i, j)} \to \gamma$ (where we use the lexical order on the $(i, j)$),
and $$\bigoplus_{ij} u_{ij} = \bigoplus_i (\bigoplus_j u_{ij}) \ge \bigoplus_i v_i \ge w.$$
Therefore $X \doms Z$. 
\end{proof}

\subsection{Two additional inequalities for domination}
For $X = \sum w_i \alpha_i \in \WW(U)$, and $A$ a component of $\di U$, let 
$$X|_A = \sum_i (\# \di \alpha_i \cap \di A) w_i,$$
where of course,
for each $i$,
$\# \di \alpha_i \cap \di AU \in 0, 1, 2$. 
\begin{lemma} \label{lem:dom-restr}
Suppose $X \doms Y$, where $X \in \WW(U)$ and $Y \in \WW(V)$, let $A$ be a component of $\di U$,
and assume $A$ is a component of $\di  V$.
Then $X|_A \ge Y|_A$. 
\end{lemma}
\begin{proof}
We can reduce to the case where $X = \sum_i  w_i \alpha_i$, $Y = v \beta$, $(\alpha_i) \to \beta$, and $\oplus_i w_i \ge v$. 
If $(\# \di \beta \cap \di A) = 1$, then there must be some $\alpha_i$ with $(\# \di \alpha_i \cap \di A) = 1$ and, moreover, $w_i \ge v$. The case where $(\# \di \beta \cap \di A) = 2$ is similar. 
\end{proof}

\begin{lemma} \label{lem:in-each}
Suppose that $X \doms Y$,
and $\beta \in \supp Y$.
Suppose that there are disjoint subsets $E_1, \ldots, E_n$ of $\supp X$,
such that whenever $(\alpha_j) \to \beta$, 
then there is an $\alpha_j$ in each of the $E_k$.
Then 
$$
Y(\beta) \le \bigoplus_k \norm{X|_{E_k}}.
$$
\end{lemma}
\begin{proof}
We have $Y(\beta) = \sum w_i$ and $X \ge \sum X_i$,  
where for each $i$, 
$X_i =  \sum v_{ij} \alpha_{ij}$,
$(\alpha_{ij})_j \to \beta$, and $\oplus_j v_{ij} \ge w_i$. 
Then,
for each $i$,
we have
$$
w_i \le \bigoplus_k \norm{X_i|_{E_k}}.
$$
Therefore, applying \eqref{eq:sum-hsum},
$$
Y(\beta) \le \sum_i \bigoplus_k \norm{X|_{E_k}} \le \bigoplus_k \sum_i \norm{X_i|_{E_k}} \le \bigoplus_k \norm{X|_{E_k}}. 
\qedhere
$$
\end{proof}

\subsection{Addition and subtraction and domination}
\begin{lem} \label{minus b}
If
$$\bigoplus ( x_i + b_i ) \ge y,$$
then
$$\bigoplus x_i \ge y - \sum b_i.$$
\end{lem}

\begin{proof}
$$
  \frac {\partial}{\partial x_i}\bigoplus x_k = \frac{(\bigoplus x_k)^2}{x_i^2} \le 1.
$$
Therefore
$$
\bigoplus (x_i + b_i) \le \bigoplus x_i \; + \; \sum b_i. \qedhere
$$
\end{proof}

\begin{lem} \label{lem:subtract-from-right}
If $X + B \lollypop Y$,
then $X \lollypop Y - \norm{B_1}_1$.
\end{lem}

\begin{proof}
Now suppose $X + B \lollypop Y$.
Formally we can write
$X + B \ge \sum_i T_i$, $Y = \sum Y_i$,
where
$T_i = \sum_j w_{ij} \alpha_{ij}$,
$Y_i = v_i \beta_i$,
and
\begin{equation} \label{arrow}
(\alpha_{ij})_j \rightarrow \beta_i,
\end{equation}
and
$$\bigoplus_j w_{ij} > v_i.$$
By the general theory of positive vectors in $\R^n$,
we can write
$T_i = X_i + B_i$,
where $X \ge \sum X_i$,
and $B \ge \sum B_i$.
So writing
$X_i = \sum_j w^X_{ij} \alpha_{ij}$,
and
$B_i = \sum_j w^B_{ij} \alpha_{ij}$,
we obtain
$$\bigoplus_j (w^X_{ij} + w^B_{ij}) \ge v_i,$$
so, by Lemma \ref{minus b},
$$\bigoplus_j w^X_{ij} \ge v_i - \norm{B_i}_1,$$
and therefore (using (\ref{arrow}))
$$
X \ge \sum_i X_i \;\lollypop\; \sum_i (v_i - \norm{B_i}_1)\beta_i \;\ge\; Y - \norm{B}_1. \qedhere
$$
\end{proof}
We then have the following corollary to Lemma \ref{lem:subtract-from-right}.
\begin{cor} \label{cor:sub-sub}
Suppose $U \subset V$, $m \in \R^+$, $X \in \WW(U)$,  $Y \in \WW(V)$,
and 
$$
X  \multimap Y.
$$
Then 
$$
X -m \multimap Y- 3 |\chi(U)| m.
$$
\end{cor}
\begin{proof}
We have some WAD $B$ on $U$ such that $X = (X - m) + B$, and $\norm{B}_\infty \le m$ (and hence $\norm{B}_1 \le 3 |\chi(U)| m$).   
\end{proof}

\renewcommand{\S}{\mathbf S}

\section{Appendix C: Estimates in hyperbolic geometry}

In this appendix we will describe how to measure the geometry of a surface
 using transverse geodesic arcs.
We will then show how to compute the lengths of peripheral closed geodesics of the surface
 using these measurements.

Suppose $T$ is a compact hyperbolic surface with geodesic boundary.
The following lemma appears as the Corollary to section 3.3 of \cite{Ab}:

\begin{lem} \label{margulis}
There is an $\ezero > 0$ such that any two distinct closed geodesics on $T$ of length at most $\ezero$
 are simple and disjoint.
\end{lem}

Let $\S$ be a compact hyperbolic surface with geodesic boundary.
Then a \emph{transverse geodesic arc} for $\S$
 is a proper path of minimal length in its proper homotopy class.
If $\bfalpha$ is a path on $\S$,
it is a transverse geodesic arc if and only if
 it is a geodesic arc that meets $\partial \S$ orthogonally,
or,
equivalently,
the double of $\bfalpha \cup \overline \bfalpha$
 in $\S \cup \overline \S$
is a closed geodesic.

Let $S$ be a compact Riemann surface with boundary,
and endow $\Int S$ with its Poincar\'e metric.
The peripheral geodesics on $\Int S$ bound a compact surface
$\S$ with geodesic boundary,
called the \emph{convex core} of $S$.
There is a homeomorphism $h\from \S \to S$
that is isotopic through embeddings to the inclusion $\S \subset S$.
We can then form a weighted arc diagram $M^S$ on $S$ as follows:
 for $\alpha \in \AA(S)$,
we find the transverse geodesic arc $\bfalpha$ for $\S$
such that $h(\bfalpha) \sim \alpha$.
Then we let $M^S(\alpha) = - \log L(\bfalpha)$ if
 $L(\bfalpha) < \ezero/2$,
and $M^S(\alpha) = 0$ otherwise.
Then $M^S$ is supported on a set of disjoint arcs,
so $M^S$ is a weighted arc-diagram for $S$.
We let
$\M^S = \bigcup\{ \bfalpha : L(\bfalpha) < \ezero/2 \}$,
so $\M^S \subset  \S$ is the union of the short transverse geodesic arcs of $\S$.

We will call a non-peripheral simple closed geodesic
 a \emph{dividing} geodesic.
The following result appears in \cite{Ab}:

\begin{lem} \label{abikoff}
Let $T$ be a compact hyperbolic Riemann surface with bounded-length geodesic boundary.
Then either $T$ is a pair of pants,
or there is a bounded-length dividing geodesic on $T$.
\end{lem}

We say that a hyperbolic surface $T$ is \emph{symmetric}
 if it admits an isometric orientation-reversing involution,
which we will denote by complex conjugation: ``$z \mapsto \overline z$''.
Then we let $E_T = \{ z \in T : z = \overline z \}$,
 and $E_T$ will be a union of (simple) closed geodesics and transverse geodesic arcs.
($E_T$ depends implicitly on on the choice of involution.
Whenever we say ``symmetric hyperbolic surface'' we will mean
 ``symmetric hyperbolic surface and choice of involution.'')
Note that $T \sm E_T$ has two components,
call them $A_T$ and $\overline{A_T}$,
which are mapped to each other by $z \mapsto \overline z$.
We prove a symmetric version of Lemma \ref{abikoff}:

\begin{lem} \label{symmetric abikoff}
Every symmetric compact hyperbolic surface $T$ with \break bounded-length geodesic boundary
 has a bounded-length symmetric pair of pants decomposition.
\end{lem}
\begin{proof}
It suffices to find a single bounded-length symmetric dividing geodesic on the surface,
or a symmetric pair of disjoint bounded-length dividing geodesics,
because then we can cut the surface $T$ along that geodesic or pair of geodesics,
and repeat.

By Theorem \ref{abikoff},
unless $T$ is a pair of pants,
there is a dividing geodesic $\gamma$ of bounded length on $T$.
If $\gamma \cap E_T = \emptyset$,
then $\gamma \cap \overline \gamma = \emptyset$,
and we are done.
Likewise, if $\gamma \subset E_T$,
then $\gamma$ is symmetric,
and we are done.
Otherwise,
let $\eta$ be a component of $\gamma \cap \Cl A_T$.
Then $\eta \cup \overline \eta$
is a non-trivial non-peripheral simple closed curve
so we let $\tau$ be the dividing
 geodesic homotopic to $\eta$.
Then $L(\tau) \le 2 L(\eta) < 2 L(\gamma)$,
so $\tau$ is the desired object.
\end{proof}

We now prove two basic estimates on transverse geodesic arcs on pairs of pants.
We denote by
$[x, y, z, r, s, t]$
the right-angled hyperbolic hexagon
with lengths $x, y, z, r, s, t$ in that order.
We will omit lengths that are not specified,
so for example $[a,,b,,c,]$
denotes the right-angled hexagon with alternating side lengths $a$, $b$, and $c$.
We first estimate the length of one side in a hyperbolic right-angled hexagon,
in terms of the lengths of three alternating sides:
\begin{lem}\label{fenchel}
Let $[a, c', b, , c, ,]$ be a hyperbolic right-angled hexagon,
and suppose that $a, b, c \le r$.
Then $c' = -\log a - \log b + O(1;r)$.\footnote{Notation $O(1;r)$ stands for a quantity  bounded
     in terms of  $r$. }
\end{lem}
\begin{proof}
We use the formula (from \cite{Fenchel:Elementary}):
\begin{equation}\label{hexagon}
\cosh c' = \frac{\cosh c + \cosh a \cosh b}{\sinh a \sinh b}
         = e^{O(1; r)} \frac 1 {ab}
\end{equation}
and recall that $\cosh^{-1} x = \log x + O(1; t)$
whenever $x \ge t > 0$.
\end{proof}
We let $\PP(a, b, c)$ denote the pair of pants with cuff lengths $a$, $b$, and $c$.
\begin{lemma}
We have the following two estimates:
\begin{enumerate}
\item
If $P = \PP(a, b, c)$ is a pair of pants,
and $\gamma$ is the transverse geodesic arc that connects $a$ and $b$,
then
$$
|\gamma| = -\log a - \log b + O(1),
$$
for $a, b, c \le C_0$.
\item
If $\gamma$ is the transverse geodesic arc that connects $a$ and $a$ in $P(a, b, b)$,
then
$$
|\gamma| = - 2 \log a + O(1)
$$
for $a, b \le C_0$.
\end{enumerate}
\end{lemma}
\begin{proof}
We prove each of the above:
\begin{enumerate}
\item
We cut $P$ along the three pairwise transversals
into two right-angled hexagons.
By formula (\ref{hexagon}) these hexagons are equal,
and hence each has type
$[a/2, \gamma,b/2, \, , c/2, \,]$.
Apply now Lemma \ref{fenchel}.

\item
We let $\eta$ be the transversal from the length $a$ cuff to one other.
Then we cut along $\eta$ and $\gamma$ to obtain the right-angled hexagon
$[a/4, \gamma, a/4, \eta, b, \eta]$,
and then apply Lemma \ref{fenchel}. \qedhere
\end{enumerate}
\end{proof}

Given a closed geodesic $\gamma$ and an arc $\alpha\in \AAA$,
let
$\pairing{ \gamma, \alpha }$
stand for the intersection number of $\gamma$ with $\alpha$,
i.e.,
the minimal number of intersections of $\gamma$ with the paths representing $\alpha$.
Given a weighted arc diagram $W=\sum W(\alpha)\alpha$,
we can define the intersection number
$\pairing{ \gamma, W }=\sum W(\alpha)\pairing{ \gamma, \alpha }$
by linearity.
We can now prove the main theorem of this appendix:
\begin{thm} \label{peripheral lengths}
Let $S$ be a compact Riemann surface with boundary,
and endow $\Int S$ with its Poincar\'e metric.
Suppose that $\gamma$ is a peripheral closed geodesic for $S$.
Then
$$
L(\gamma) = 2 \left< M^S, \gamma \right> + O(1; \chi(S)).
$$
\end{thm}

\begin{proof}
We let $\S$ be the convex core of $S$.
We find a symmetric bounded-length pair of pants decomposition for $\S \cup \overline \S$
extending $\M^S\cup \bar \M^S$.
Then we can write $\gamma = \bigcup t_i$,
where the segments $t_i$ are interior-disjoint,
and each is a transverse arc of one of the pairs of pants.
Then $L(t_i) = -\log a_i - \log b_i + O(1)$,
where $a_i$ and $b_i$ are the lengths of the cuffs that $\gamma$ connects
(possibly the same cuff).
Therefore
$$
L(\gamma) = 2 \sum - \log a_i + O(1)
$$
where the $a_i$ are the lengths of the cuffs that $\gamma$ crosses,
counted with multiplicity.
But
\begin{equation}
2 \sum -\log a_i = 2 \left< M^S, \gamma \right> + O(1). \qedhere
\end{equation}
\end{proof}

\end{document}